\documentclass{article}
\usepackage{authblk}
\usepackage[utf8]{inputenc}
\usepackage{url}
\usepackage[version=4]{mhchem}
\usepackage{multirow}
\usepackage{lipsum}
\usepackage{stackrel}
\usepackage{listings}
\usepackage{amssymb}
\usepackage{nccmath}
\usepackage{a4wide}
\usepackage{mathtools}
\usepackage{amsmath,amsfonts,amsthm,amssymb}
\usepackage{bbm, dsfont}
\usepackage[justification=justified, singlelinecheck=false]{caption}
\usepackage{graphicx}
\usepackage{enumitem}
\usepackage[justification=centering]{caption}
\usepackage{algorithmicx}
\usepackage{enumitem}
\usepackage{empheq,etoolbox}
\patchcmd{\subequations}
  {\theparentequation\alph{equation}}
  {\theparentequation.\alph{equation}}
  {}{}
\usepackage[justification=centering]{caption}
\usepackage{algorithmicx}
\usepackage{algorithm}
\usepackage{algpseudocode}
\usepackage{stackengine}

\usepackage[title]{appendix}
\usepackage[scaled=.90]{helvet}
\usepackage{courier}
\usepackage{ae}
\usepackage[T1]{fontenc}
\usepackage{graphicx}
\usepackage{subcaption}
\usepackage[export]{adjustbox}
\usepackage{wrapfig}

\usepackage[T1]{fontenc}
\usepackage{here} 
\usepackage[colorlinks=false, pdfborder={0 0 0}]{hyperref}

\usepackage{mdwlist}
\usepackage{tcolorbox}
\usepackage{fancybox}
\usepackage{framed}

\definecolor{darkblue}{rgb}{0,0,0.8}
\definecolor{darkgreen}{rgb}{0,0.8,0}

\definecolor{magenta}{rgb}{0.5,0,0.5}

\newcommand{\mathleft}{\@fleqntrue\@mathmargin0pt}

\newcommand{\LZ}{\vspace{\baselineskip}}

\newcommand{\HLZ}{\vspace{0.5\baselineskip}}

\newtheorem{remark}{Remark}[section]

\providecommand{\keywords}[1]
{
  \small	
  \textbf{\textit{Keywords---}} #1
}

\begin{document}
\title{Exact simulation of the first-passage time of SDEs to time-dependent thresholds} 
\author[1]{Devika Khurana\thanks{Email: devika.khurana@jku.at}}
\author[2]{Sascha Desmettre\thanks{Email: sascha.desmettre@jku.at}}
\author[1]{Evelyn Buckwar\thanks{Email: Evelyn.Buckwar@jku.at}}
\affil[1]{\centerline{\small Institute of Stochastics, Johannes Kepler University Linz} }
\affil[2]{\centerline{\small Institute of Financial Mathematics and Applied Number Theory, Johannes Kepler University Linz} }
\date{\today}

\maketitle
\begin{abstract}
 \noindent The first-passage time (FPT) is a fundamental concept in stochastic processes, representing the time it takes for a process to reach a specified threshold for the first time. Often, considering a time-dependent threshold is essential for accurately modeling stochastic processes, as it provides a more accurate and adaptable framework. In this paper, we extend an existing Exact simulation method developed for constant thresholds to handle time-dependent thresholds. Our proposed approach utilizes the FPT of Brownian motion and accepts it for the FPT of a given process with some probability, which is determined using Girsanov's transformation. This method eliminates the need to simulate entire paths over specific time intervals, avoids time-discretization errors, and directly simulates the first-passage time. We present results demonstrating the method's effectiveness, including the extension to time-dependent thresholds, an analysis of its time complexity, comparisons with existing methods through numerical examples, and its application to predicting spike times in a neuron.\HLZ\\
\keywords{First-passage time, \and Exact simulation, \and Stochastic differential equations, \and Neuroscience} \textbf{\textit{MSC 2010 Classification}}---- 37M05, 65C30, 60G05, 60H35, 68Q87, 92C20
\end{abstract} 


\section{Introduction}
Dynamical systems are frequently modeled using one-dimensional diffusion processes in various fields. Understanding their behavior involves examining the time it takes for the diffusion to reach a specified threshold. The time at which the diffusion reaches the threshold for the first time is known as the first-passage time.
For instance, in financial markets, the pricing of barrier options, cf. \cite{barone2008barrier}, depends on whether the system crosses the threshold before the option's expiry time.
In oncology, first-passage time analysis is critical for assessing treatment efficacy. A notable application is the study of tumor growth delay \cite{roman2021using}, where first-passage time metrics are used to compare the time required for a tumor to reach a certain volume (threshold) with and without the treatment.

Additionally, the theory of first-passage time extends to diverse domains, including physics \cite{PhysRevE.62.6065}, reliability engineering \cite{chun2019system}, and neuroscience. In our focus on the neuroscience model \cite{kostal2008efficient}, the diffusion process describes the evolution of the membrane voltage of a neuron. When the membrane voltage reaches a certain threshold, it generates a spike, and the information is transmitted to the next neuron. The time at which the membrane voltage hits the threshold is known as the spike time, which can be modeled as the first-passage time.\HLZ\\
The corresponding mathematical problem can be formulated as follows:\\
Let $\left(\Omega,\mathcal{F},\mathbb{P}\right)$ be a probability space with filtration $\mathbb{F}=\{\mathcal{F}_{t},t\geq0\}$. We consider the stochastic process $\{Y_{t},t \geq 0\}$ adapted to the filtration $\mathbb{F}$ as a solution to a time-homogeneous It{\^o} stochastic differential equation: 
\begin{equation}
     dY_{t}=\mu\left(Y_{t}\right)dt+\sigma\left(Y_{t}\right)dB_t, \hspace{1cm}   Y_{0}=y_0\in \mathbb{R}\,,
     \label{eq:e1}
\end{equation}
where $\{B_{t},t\geq0\}$ is the standard one-dimensional Brownian motion adapted to $\mathbb{F}$.\HLZ\\
\textbf{Assumption 1.1} The drift $\mu$ and diffusion $\sigma$ of SDE described in (\(\ref{eq:e1}\)) are defined in a manner that ensures the existence of a unique strong solution.\HLZ\\
We denote $\tau_{\beta}$ as the first-passage time of the diffusion process $Y$ to a time-dependent deterministic threshold $\beta\left(t\right)$, where $\beta:\mathbb{R}^{+}_{0}\rightarrow\mathbb{R}$. The first-passage time $\tau_\beta$, inherits randomness from $Y$ and thus is a random variable defined as follows:
\begin{equation}
    \tau_\beta:=\inf\{t\geq0:Y_t=\beta(t)\}.
    \label{eq:e2}
\end{equation}
Due to randomness inherited by $\tau_\beta$, path dependence in the process $Y$, and the dynamic nature of the threshold function over time, determining the first-passage time can be tedious.\\
The existing methods for evaluating the first-passage time $\tau_\beta$, based on the dynamics in which this time is obtained, can be classified as follows:

$\left(i\right)$\textit{ Explicit expressions of probability density function} for the first-passage time of a stochastic differential equation to a time-dependent threshold function, or even to a constant one, are exceptionally rare. When available, these expressions often correspond to basic mathematical models that often lack direct applicability \cite{jeanblanc2009mathematical}.

$\left(ii\right)$\textit{ Approximation of the probability density function} of the first-passage time. One prominent approach is to use the Fokker-Planck equations \cite{oksendal2013stochastic}. The probability that the process $Y$ reaches the threshold $\beta\left(t\right)$ before $T$ is given by $\phi_{y}=F\left(y,T\right)$, where the function $F\left(y,T\right)$ satisfies a Fokker-Planck equation. Given the complexity of such equations, exact solutions are rarely feasible. A novel strategy as explored in \cite{boehm2021fast} involved transforming the space-time domain to $\left(0,1\right)\times\left(0,1\right)$ and using a fast-converging series expansion for the first-passage time distribution to simplify the Fokker-Planck equation. More work in this direction can be found in \cite{boehm2022efficient,rasanan2023numerical}. Another noteworthy approach focuses on finding the transitional density $g(\beta(\tau_{\beta}),\tau_{\beta}|y_{0},t_{0})$ of the FPT $\tau_{\beta}$. This density is the solution to a Volterra-type integral equation \cite{durbin1971boundary,van1992stochastic}. However, the kernel of this equation is weakly singular, making numerical approximation inherently unstable \cite{durbin1971boundary}. Additionally, one can only compute the solution of the equation for one initial value of the diffusion process at a time. Relevant studies on this topic include \cite{giorno2021first,valov2009integral,ricciardi1984integral}. An alternative approach is the Method of Images \cite{daniels1982sequential}. This method involves the construction of a measure $A$ using a linear combination of simpler measures. This construction facilitates an approximate threshold for which the first-passage time distribution can be computed \cite{zipkin2016method}. Another approach is transforming the pair of SDE and threshold to another pair for which the pdf of first-passage time can be easily computed \cite{cherkasov1957transformation,capocelli1976transformation,ccelebiouglu2013transformation}.

$\left(iii\right)$\textit{Approximation of the value of the first-passage time} can be obtained by simulating the path of the process $Y$ and then determining the time at which the simulated path crosses $\beta\left(t\right)$. The path of $Y$ can be simulated using various time-discretisation schemes. The most classical one is the Euler-Maruyama method, which is a widely used approach due to its simplicity in implementation. However, it tends to overestimate the first-passage time if the step size is not sufficiently small \cite{buchmann2003computing}. Furthermore, reducing the step size to decrease error leads to increased variance. That is why there is always a challenge to balance the time-discretization and estimation errors. One improvement was introduced in \cite{mannella1999absorbing}, where the probability that the process $Y$ crosses the threshold within a time step was established, and later Gobet in \cite{PS_2001__5__261_0} proved that the method has a $o\left(1/2\right)$ time discretization error, compared to the $o\left(1\right)$ error of the classical Euler-Maruyama method. Other notable improvements can be found in \cite{ichiba2011efficient,gobet2010stopped,giraudo2001monte}. Another technique, the random walk on spheres introduced by Muller \cite{muller1956some}, focuses on a multi-dimensional Brownian motion $Z$. It aims to determine the first exit time of $Z$ from its domain $D$. This method begins by constructing the largest sphere centered at initial value $Z_{0}$ contained within $D$. A point is then chosen uniformly from the boundary of this sphere to serve as the hitting position of the sphere, denoted $z_{1}$. Next, a new sphere is constructed, centered at $z_{1}$ and contained within $D$. The exit time from this new sphere is determined by uniformly selecting a point from its boundary. This process is repeated iteratively, with each new sphere being the largest possible one centered at the current point. The algorithm continues until the exit position is sufficiently close to the boundary $\partial D$. More work in this direction can be found in \cite{deaconu2006random,deaconu2017walk}.

$\left(iv\right)$\textit{Exact simulation}, as introduced by Beskos and Roberts in \cite{beskos2005exact} to simulate the exact grid of the process $Y$, uses the acceptance-rejection method. The fundamental idea is to find another process that closely resembles $Y$ and is easier to simulate. A finite number of points from the path of this auxiliary process are then accepted or rejected as points on the path of $Y$ based on a probability derived using Girsanov's transformation \cite{jeanblanc2009mathematical}. The remaining path between the accepted points is then constructed using the Brownian Bridge. To obtain the first-passage time, one can determine the point at which this path crosses the threshold. However, Hermann and Zucca in \cite{herrmann2019exact} introduced a technique that directly and exactly simulates the first-passage time of $Y$ to a constant threshold. Similar to the previous approach, this method begins by identifying a process close to $Y$. The sample of the first-passage time of this auxiliary process to $\beta$ (constant threshold) is then accepted as the first-passage time of $Y$ to $\beta$, with a probability obtained using Girsanov's transformation. This approach is particularly appealing because it eliminates the need to simulate the entire path and incurs no approximation error but only inherent statistical error. Moreover, it does not restrict the first-passage time to be within a fixed time interval, which is advantageous given that first-passage times can often be very large - an issue common to other methods. More work in this direction can be found in \cite{herrmann2020exact,herrmann2023exact,beskos2006retrospective,blanchet2020exact}. \HLZ\\
While the Exact simulation method is notable due to the absence of time-discretization errors, it traditionally applies only to scenarios with a constant threshold. Recognizing this constraint is crucial, as many practical models operate with time-varying thresholds, cf. \cite{sacerdote2005inverse,gutierrez2013american,tamborrino2016approximation,levakova2019adaptive}.\HLZ\\
This paper extends the Exact simulation method by making it applicable to scenarios where the threshold function is time-dependent. We use the first-passage time of Brownian motion as the auxiliary process to implement acceptance-rejection sampling to obtain the first-passage time of the process $Y$. For certain threshold functions, such as linear and piecewise linear functions, we have explicit density functions from which to sample. For other types of curvy threshold functions, we employ the algorithmic approach introduced in \cite{herrmann2016first} to approximate the first-passage time of Brownian motion. This approach is efficient in terms of convergence rate and computational cost.\HLZ\\
While the first-passage time can be perceived as a time value, it is essential to recognize that this value varies for different potential paths of the stochastic differential equation, introducing an element of randomness. We need the following standing assumption: 
\HLZ\\
\textbf{Assumption 1.2} The first-passage time $\tau_{\beta}$ is almost surely finite and its probability density function $f_{\tau_\beta}$ is well-defined.\HLZ\\
The structure of the paper is organized as follows: Section 2 introduces the rejection sampling method; in particular, it discusses the calculation of the rejection probability and the corresponding algorithm for accepting or rejecting samples of the first-passage time of Brownian motion for $\tau_{\beta}$. Section 3 applies the methodology to a linear threshold and then to a curved threshold with the help of examples. In Section 4, we analyze the time complexity of the algorithm and the scope of reducing it. Finally, Section 5 applies the method to predict spike times in a neuron, where the membrane voltage is modeled using a quadratic leaky integrate-and-fire model with an exponentially declining adaptive threshold.

\section{Rejection Sampling for the FPT to Moving Boundaries}
In this section, we formulate the Exact Algorithm for simulating $\tau_{\beta}$ using rejection sampling. The core idea is to identify a process close to $\tau_{\beta}$ (we choose the first-passage time of Brownian motion), extract a sample from this new process, and then determine the probability of rejecting this sample as a valid realization for $\tau_{\beta}$.\\ 
To compute the probability of rejection, we use Girsanov's transformation, specifically the version dealing with stopping times.\HLZ\\
As this version of Girsanov's transformation applies to processes satisfying a stochastic differential equation (SDE) with unit diffusion, we focus on the case when the diffusion function of the SDE is one or can be transformed to one. We define the process $X_t:=F(Y_t)$, where
\[
F\left(Y_{t}\right):=\int_{}^{Y_{t}}\frac{dz}{\sigma\left(z\right)}
\]
and $Y_t$ satisfies equation (\ref{eq:e1}).\HLZ\\
\textbf{Assumption 2.1} The diffusion coefficient $\sigma$ in equation (\ref{eq:e1}) is non-negative and the function $F$ is bijective.\HLZ\\
With Assumption 2.1, we can use Lamperti's transform as defined in~\cite{panik2017stochastic} and write $\{X_{t},t\geq 0\}$ as the solution to the transformed SDE, 
\begin{equation}
    dX_{t}=\alpha\left(X_{t}\right)dt+dB_{t}, \hspace{0.5cm} X_{0}=x_0,
    \label{eq:e3}
\end{equation}
where the transformed drift function $\alpha\left(X_{t}\right)$ has the form
\[
\alpha\left(x\right)=\frac{\mu\left(F^{-1}\left(x\right)\right)}{\sigma\left(F^{-1}\left(x\right)\right)}-\frac{1}{2}\sigma_{X}\left(F^{-1}\left(x\right)\right).
\]
Both $X$ and $B$ are defined on the same probability space $(\Omega,\mathcal{F},\mathbb{P})$ as $Y$.\\
Henceforth, we will consider (\ref{eq:e3}) as the underlying SDE from now on.\HLZ\\
\textbf{Theorem 2.1} \textbf{(Girsanov Transformation for Stopping Times)} Let $\rho$ be the explosion time of the solution of the SDE described in (\ref{eq:e3}) such that $\lim_{t\to\rho}X_t=\infty$. Define a measure $\mathbb{Q}$ by
\[
d\mathbb{Q}(\omega) = M_{\mathcal{T}}(\omega)d\mathbb{P}(\omega),
\]where \[
M_\mathcal{T} = \exp{\left(\int_{0}^{\mathcal{T}}\alpha(\omega)dB_s-\frac{1}{2}\int_{0}^{\mathcal{T}} \alpha^{2}(\omega)ds\right)},
\]
and $\mathcal{T}<\rho$ is any stopping time.\HLZ\\
Assume that $M_t$ is a martingale on $0\leq t \leq \mathcal{T}$ w.r.t. $\mathbb{P}$. Then, under $\mathbb{Q}$, $X$ is a Brownian motion for $0\leq t \leq \mathcal{T}$.\HLZ\\
Consequently, there exists a probability measure $\mathbb{Q}$, under which $X$ is a Brownian motion such that for any bounded measurable function $\Psi$,
\[
\mathbb{E}_{\mathbb{P}}[\Psi(t \leq\mathcal{T})]=\mathbb{E}_{\mathbb{Q}}[\Psi(t \leq\mathcal{T})M_{\mathcal{T}}].
\]
For a detailed description and proof of this theorem, please refer to \cite{jeanblanc2009mathematical,oksendal2013stochastic} and the references therein.\qed
\begin{remark}
    The Novikov condition
\begin{equation}
    \mathbb{E}_{\mathbb{P}}\left[\exp\left(\frac{1}{2}\int_{0}^{\mathcal{T}}\alpha^{2}(\omega)ds\right)\right]<\infty
    \label{eq:e4}
\end{equation}
is sufficient to ensure that $\{M_t,0\leq t\leq \mathcal{T}\}$ is a martingale.
\end{remark}
Using Theorem 2.1, we show that the FPT of the process under one measure has the same law as the FPT of Brownian motion under a different measure, with some probability weight. This transformation facilitates a link between the FPT of our chosen process (which is a Brownian motion) and the desired process $\tau_{\beta}$.

\subsection{Rejection Probability using Girsanov's Transformation}

In this subsection, we aim to establish a link between the first-passage time of the process $X$ and the first-passage time of a Brownian motion defined w.r.t. the measure $\mathbb{Q}$ as given by Theorem 2.1 to the threshold $\beta(t)$. Consequently, we will determine the probability weight to accept or reject a sample of the latter for the former.\HLZ\\
Let us from now on denote $X^{\mathbb{P}}$ as the process defined on $(\Omega, \mathcal{F}, \mathbb{P})$ satisfying SDE described in equation (\ref{eq:e3}) and $X^{\mathbb{Q}}$ is the same process defined on $(\Omega, \mathcal{F}, \mathbb{Q})$ where $\mathbb{Q}$ is the measure described in Theorem 2.1. Note that by Theorem 2.1 $X^{\mathbb{Q}}$ is a Brownian motion on $[0,\tau_{\beta}]$, supposed that the density process $M_{\tau_{\beta}}$ is a martingale.\HLZ\\
We need the following additional assumptions:\HLZ\\
\textbf{Assumption 2.2} (A1) The drift coefficient $\alpha(X_t)$ is square integrable.

This assumption ensures that the It\^{o} integral of $\alpha(X_t)$ is well-defined.\HLZ\\
(A2) The drift $\alpha:\mathbb{R}\rightarrow \mathbb{R}$ is a $C^{1}$ function.\HLZ\\
(A3) The threshold $\beta: \mathbb{R}^{+}_{0} \rightarrow \mathbb{R}$ is differentiable. \HLZ\\ 
Let us now define the following functions:
\begin{align*}
    &\gamma_1: [0, \infty) \rightarrow \mathbb{R} & &\gamma_2: \mathbb{R} \rightarrow \mathbb{R} \\
    &\gamma_1(t) := -\alpha(\beta(t))\beta'(t), & &\gamma_2(x) := \frac{1}{2}(\alpha'(x) + \alpha^{2}(x)).
\end{align*}
The next theorem links the first-passage time of $X^{\mathbb{P}}$ and $X^{\mathbb{Q}}$, and provides the rejection probability for sampling.\LZ\\
\textbf{Theorem 2.2} If Assumption 2.2 holds, then for any bounded measurable function $\Psi: \mathbb{R}\rightarrow\mathbb{R}$, we have that
     \begin{equation}
         \mathbb{E}_{\mathbb{P}}\left[\Psi \left(\tau_{\beta}\right)\mathbf{1}_{\{\tau_{\beta}<\infty\}}\right]= e^{-A\left(x\right)+A\left(\beta\left(0\right)\right)}\mathbb{E}_\mathbb{Q}\left[\Psi \left(\tau_{\beta}\right) \eta\left(\tau_{\beta}\right)\right],
         \label{eq:e5}
     \end{equation}
where \begin{equation}
        \eta\left(t\right):=\mathbb{E}_\mathbb{Q}\left[\exp\left(-\int_0^{t} \left(\gamma_{1}\left(s\right)+\gamma_{2}\left(X^{\mathbb{Q}}_{s}\right)\right)ds\right)|\tau_{\beta}=t\right].
        \label{eq:e6}
      \end{equation}   
Here, $\eta(t)$ represents the probability weight for rejection sampling.\LZ\\    
\textit{Proof} 
To prove (\ref{eq:e5}), we use Girsanov's transformation described in Theorem 2.1, resulting in the following equality:
      \begin{equation*}
          \mathbb{E}_{\mathbb{P}}\left[\Psi \left(\tau_{\beta}\right)\mathbf{1}_{\{\tau_{\beta}<\infty\}}\right]= \mathbb{E}_\mathbb{Q}\left[\Psi \left(\tau_{\beta}\right)\exp\left(\int_0^{\tau_{\beta}}\alpha\left(X^{\mathbb{Q}}_{s}\right)dB_{s}-\frac{1}{2} \int_0^{\tau_{\beta}} \alpha^{2}\left(X^{\mathbb{Q}}_{s}\right)ds\right)\right].
      \end{equation*}
Under Assumption 2.2 (A2), there exists a function $A$ such that $A'=\alpha$. Additionally, applying Assumption 2.2 (A1), we can use the It\^{o} formula on $A$, yielding
\[
\displaystyle{\int_0^{\tau_{\beta}}}\alpha\left(X^{\mathbb{Q}}_{s}\right)dB_{s}=A\left(X^{\mathbb{Q}}_{s}\right)_{0}^{\tau_{\beta}}-\frac{1}{2}\int_0^{\tau_{\beta}}\alpha'\left(X^{\mathbb{Q}}_{s}\right)ds,
\] which allows us to rewrite the equation as:
\begin{equation*}
\begin{split}
    \mathbb{E}_{\mathbb{P}}\left[\Psi \left(\tau_{\beta}\right)\mathbf{1}_{\{\tau_{\beta}<\infty\}}\right] = \mathbb{E}_\mathbb{Q}\Biggl[\Psi \left(\tau_{\beta}\right)\exp\Biggl(&A(X^{\mathbb{Q}}_{\tau_{\beta}})-A(X^{\mathbb{Q}}_{0})-\\
    &\frac{1}{2}\int_0^{\tau_{\beta}}\alpha'\left(X^{\mathbb{Q}}_{s}\right)ds-\frac{1}{2} \int_0^{\tau_{\beta}} \alpha^{2}\left(X^{\mathbb{Q}}_{s}\right)ds\Biggr)\Biggr]\,.
\end{split}
\end{equation*}
Under the measure $\mathbb{Q}$, $\tau_{\beta}=\inf\{t\geq 0:X^{\mathbb{Q}}_{t}=\beta(t)\}$ , which means $X^{\mathbb{Q}}$ reaches $\beta\left(t\right)$ at $\tau_{\beta}$ for the first time. Therefore, $X^{\mathbb{Q}}_{\tau_{\beta}}=\beta(\tau_{\beta})$, and we obtain\HLZ\\
$\mathbb{E}_{\mathbb{P}}\left[\Psi \left(\tau_{\beta}\right)\mathbf{1}_{\{\tau_{\beta}<\infty\}}\right]$
\begin{align*}
    & = \mathbb{E}_\mathbb{Q}\left[\Psi \left(\tau_{\beta}\right)\exp\left(A\left(\beta\left(\tau_{\beta}\right)\right)-A\left(x_{0}\right)-\frac{1}{2}\int_0^{\tau_{\beta}}\left(\alpha'\left(X^{\mathbb{Q}}_s\right)+\alpha^{2}\left(X^{\mathbb{Q}}_{s}\right)\right)ds\right)\right]\\
    & = \mathbb{E}_\mathbb{Q}\left[e^{-A\left(x_{0}\right)}\Psi \left(\tau_{\beta}\right)  \exp\left(\int_0^{\tau_{\beta}} A'\left(\beta\left(s\right)\right)\beta'\left(s\right)ds+A\left(\beta\left(0\right)\right) -\int_0^{\tau_{\beta}}\gamma_{2}\left(X^{\mathbb{Q}}_s\right)ds\right)\right]\\
    &= \mathbb{E}_{\mathbb{Q}}\left[e^{-A\left(x_{0}\right)+A\left(\beta\left(0\right)\right)}\Psi \left(\tau_{\beta}\right) \exp\left(-\int_0^{\tau_{\beta}} \left(\gamma_{1}\left(s\right)+\gamma_{2}\left(X^{\mathbb{Q}}_s\right)\right)ds\right)\right]\\
    &= \mathbb{E}_{\mathbb{Q}}\left[e^{-A\left(x_{0}\right)+A\left(\beta\left(0\right)\right)} \Psi \left(\tau_{\beta}\right) \eta\left(\tau_\beta\right)\right].
\end{align*}\qed \HLZ\\
Equation (\ref{eq:e5}) indicates that $\eta(t)$ acts as a weighting factor, adjusting the likelihood of "hitting the threshold for the first time at $\tau_\beta$" under the original measure $\mathbb{P}$ when analyzed within the framework of the new measure $\mathbb{Q}$. In cases where computing the expectation under $\mathbb{P}$ is challenging but more straightforward under $\mathbb{Q}$, $\eta(t)$ enables re-weighting of the samples, ensuring that the expectation under $\mathbb{Q}$ corresponds to the desired expectation under $\mathbb{P}$.\HLZ\\
The following section illustrates how we use this probability to perform rejection sampling.
\subsection{The Algorithm for Rejection Sampling}
Let us denote $\tau_{\beta}^{W}$ as the first-passage time of the Brownian motion $X^{\mathbb{Q}}$ to $\beta(t)$.\HLZ\\
By (\ref{eq:e5}), $\eta$ is a probability only if the following assumption holds:\HLZ\\
\textbf{Assumption 2.3} $\gamma_1\left(t\right)+\gamma_2(X^{\mathbb{Q}}_t)\geq0$ for all $0\leq t<\infty$.
\begin{remark}[Interpretation of Assumption 2.3]
    \begin{itemize}
    \item[1.] If $\gamma_1\left(t\right)\geq 0$ for $t \in \left[0,\infty \right)$ and $\gamma_2(X^{\mathbb{Q}}_t)\geq0$  for $t \in \left[0, \infty\right]$ (or $X^{\mathbb{Q}}_t \in \left(-\infty,\infty \right)$), then $\gamma_1\left(t\right)+\gamma_2(X^{\mathbb{Q}}_t)\geq0$. Thus, to proceed with the algorithm, it would be feasible to assume the non-negativity of both the functions $\gamma_1$ and $\gamma_2$, separately.
    \item[2.] This assumption needs to hold only for \(t \in [0, \tau_{\beta}^{W}]\) in the context of our algorithm.
    \item[3.] If the functions $\gamma_1$ and $\gamma_2$ also have negative values in their ranges but possess lower bounds, we can shift the functions to ensure they are non-negative. Let us assume that the functions can be negative, and their lower bounds $k_1$, $k_2$ exist such that,
\[
k_1 \leq \underset{t\in [0,\tau_\beta^{W}]}{\inf} \gamma_1 \left(t\right) \text{  and  } k_2 \leq \underset{t\in [0,\tau_\beta^{W}]}{\inf} \gamma_2 (X^{\mathbb{Q}}_t).
\]
Note that $k_1$ always exists.\\
We then define new functions as follows:
\[
\gamma_{1}^\star \left(t\right):= \gamma_1 \left(t\right)-k_1 \text{  and  } \gamma_{2}^\star (X^{\mathbb{Q}}_t):= \gamma_2 (X^{\mathbb{Q}}_t)-k_2.
\]
\end{itemize}
\end{remark}
Without loss of generality, we consider $\beta\left(0\right)>x_{0}$ and explore how we can use the rejection probability $\eta(t)$ to get the first-passage time of $X^{\mathbb{P}}$ to $\beta\left(t\right)$.\HLZ\\
The first step is to find the first-passage time of standard Brownian motion to $\beta(t)$. The explicit density distribution of the first-passage time of Brownian motion to $\beta(t)$ is available only for three cases, i.e.,  when $\beta(t)$ is 
\begin{itemize}
    \item[1)] a constant $=\beta$, the pdf is Inverse Gaussian such that $\tau_{\beta}^{W}\sim \beta/g^{2}$; $g\sim \mathcal{N}(0,1)$.
    \item[2)]linear $=at+b$, the pdf is $\sim IG(-b/a,b^2)$ if $ab<0$ and $\sim e^{-2ab}IG(b/a,b^2)$ if $ab>0$; where $IG(\lambda_{1},\lambda_{2})$ is a Inverse Gaussian distribution with parameters $\lambda_{1}$ and $\lambda_{2}$, cf. \cite{karatzas1991brownian}.
    \item[3)] piece-wise linear, cf. \cite{jin2017first}.
\end{itemize}
For other types of threshold functions, we explore a method that gives a promising approximation of $\tau_{\beta}^{W}$, which will be discussed in Section 3.1.\HLZ\\
We now present a theorem that articulates the idea of using $\eta$ as probability weight.\LZ\\
    \textbf{Theorem 2.3} Let $\tau_{\beta}^{W}$ be any positive real value. If $\Phi$ is a homogeneous Poisson process of unit intensity on the bounded region
$B=B_1\times B_2 \subset[0,\tau_{\beta}^{W}]\times\mathbb{R}^{+}_{0}$, such that $B_1=[0,\tau_{\beta}^{W}]$ and $B_2$ contains $\{\gamma_1(t)+\gamma_{2}(X^{\mathbb{Q}}_t);t\in[0,\tau_{\beta}^{W}]\}$, and $N$ is the number of points of $\Phi$ found below the graph
$\{(t,\gamma_{1}(t)+\gamma_{2}(X^{\mathbb{Q}}_{t}));t\in[0,\tau_{\beta}^{W}]\}$, then 
\[
\mathbb{P}[N=0|X^{\mathbb{Q}}_{t\leq\tau_{\beta}^{W}}]=\exp\left(-\int_{0}^{\tau_{\beta}^{W}} (\gamma_{1}(t)+\gamma_{2}(X^{\mathbb{Q}}_{t}))dt\right).
\]
\textit{Proof} We use the definition of a Poisson process and the thinning method, cf. \cite{lewis1979simulation}.\HLZ\\ Given $X^{\mathbb{Q}}_{t\leq\tau_{\beta}}$, the probability that $N$ points lie below the graph is given by \[\frac{1}{N!}\left(\int_{0}^{\tau_{\beta}^{W}}(\gamma_{1}(t)+\gamma_{2}(X^{\mathbb{Q}}_t))dt\right)^{N}\exp\left(-\int_{0}^{\tau_{\beta}^{W}}(\gamma_{1}(t)+\gamma_{2}(X^{\mathbb{Q}}_t))dt\right).\] Thus the probability that $  N=0$ is given by \[\exp\left(-\int_{0}^{\tau_{\beta}^{W}}(\gamma_{1}(t)+\gamma_{2}(X^{\mathbb{Q}}_t))dt\right).\] This implies that $N$ follows a Poisson distribution with mean \[\int_{0}^{\tau_{\beta}^{W}}(\gamma_{1}(t)+\gamma_{2}(X^{\mathbb{Q}}_t))dt.\] For further details, we refer to \cite{beskos2005exact,beskos2006retrospective}.\qed \LZ\\
Theorem 2.3 suggests that it is possible to carry out rejection sampling without evaluating $\eta(t)$.
\HLZ\\ 
This means in a simulation context: Under Assumption 2.3,  $\eta(t)$ lies in the interval $\left[0,1\right]$. By removing the exponential and the negative sign from $\eta$, the transformed $\eta(t)$ will lie within $B_2$. We thus generate a uniform random number $u$ on $B_2$ and compare it with the transformed probability. If for each time $t\in [0,\tau_{\beta}^{W}]$, $u\leq\gamma_{1}(t)+\gamma_{2}(X^{\mathbb{Q}}_{t})$, then it implies that the expected value of its negative counterpart over time is greater than $u$, indicating that the sample should be rejected. Conversely, if for each $t\in [0,\tau_{\beta}^{W}]$,  $u> \gamma_{1}(t)+\gamma_{2}(X^{\mathbb{Q}}_{t})$, we should accept the sample.\HLZ\\
To formalize this idea, we define the stochastic domain $D$ as follows:
\[
D:=\left\{(t,u)\in [0,\tau_{\beta}^{W}]\times B_2: u\leq\gamma_{1}(t)+\gamma_{2}(X^{\mathbb{Q}}_t)\right\},
\]
where $X^{\mathbb{Q}}$ is a Brownian motion with $X^{\mathbb{Q}}_{0}=x_{0}$ and $X^{\mathbb{Q}}_{\tau_{\beta}^{W}}=\beta(\tau_{\beta}^{W})$. The domain $D$ then represents the area of rejection.\\
For sampling, we need the Brownian motion trajectories to satisfy the conditions, that it starts at $x_{0}$ and ends at $\beta(\tau_{\beta}^{W})$ on $[0,\tau_{\beta}^{W}]$, where $\tau_{\beta}^{W}$ is the first-passage time. One possible approach to simulate such Brownian motion without compensating on its dynamics is to use the Bessel bridge, \[R_{s}=\beta(\tau_\beta-s)-X^{\mathbb{Q}}_{\tau_\beta-s}, \hspace{0.5cm} s\in[0,\tau_{\beta}^{W}],\] (for details see Appendix A.). Because of the relation between $X^{\mathbb{Q}}$ and $R$, the expectation $\eta\left(t\right)$ can alternatively be expressed in terms of the Bessel Bridge.\\
One could also think of using a Brownian Bridge, but a Brownian Bridge has the drawback that it can hit $\beta(t)$ before $\tau_{\beta}^{W}$, so it is not useful in our setting.\HLZ\\
Based on Theorem 2.3 and the description of the simulation described above, we present the following algorithm to simulate the first-passage time $\tau_{\beta}$ to a time-dependent threshold $\beta(t)$:
\begin{algorithm}[h]
    \caption{(Schematical Description)}
    \textbf{Input:} A sample of non-negative random variable $\tau_\beta^{W}$ using the pdf or its approximation.\vspace{0.1cm}\\
    \textbf{Step 1.} Generate a Bessel process $R$ over the time interval $[0,\tau_{\beta}^{W}]$ with endpoint\\\hspace*{1.14cm}$R_{\tau_{\beta}^{W}}=\beta\left(0\right)-x_{0}$.\vspace{0.1cm}\\
    \textbf{Step 2.} Generate a Poisson point process $\Phi$ on the space $[0,\tau_{\beta}^{W}]\times B_2$, independent of the \hspace*{1.06cm} Bessel process, with unit intensity.\vspace{0.1cm}\\
    \textbf{Step 3.} \textbf{if} $\Phi\left(D\right)=0$ \textbf{then}\vspace{0.1cm}\\
    \textbf{Step 4.} \hspace{0.36cm}assign $\tau=\tau_{\beta}^{W}$\vspace{0.1cm}\\
    \textbf{Step 5.} \textbf{otherwise} reject the current input and repeat the procedure with a new input.\vspace{0.1cm}\\
    \textbf{Output}: A random variable $\tau$.
\end{algorithm}

\textbf{Theorem 2.4} Let $f_{\tau_{\beta}^{W}}$ be the density of first-passage time of the Brownian motion $X^{\mathbb{Q}}$. If Assumptions 2.2 and 2.3 hold true, then the pdf of the outcome variable $\tau$, denoted by $f_{\tau}$ satisfies
\begin{equation}
    f_{\tau}\left(t\right)=\frac{1}{N_c} \eta\left(t\right) f_{\tau_{\beta}^{W}}\left(t\right),
    \label{eq:eq7}
\end{equation}
where $N_c$ is the normalizing coefficient
\begin{equation*}
    N_c:=\int_{0}^{\infty} \eta\left(t\right) f_{\tau_{\beta}^{W}}\left(t\right) dt.
\end{equation*}
\textit{Proof}
Using the rejection sampling described before and the concept of acceptance-rejection, we can deduce that the density of the outcome of Algorithm 1, $\tau$, is given by (\ref{eq:eq7}). By comparing with (\ref{eq:e5}), the normalizing constant, $N_c$ is $e^{A\left(x_{0}\right)-A\left(\beta\left(0\right)\right)}$.\qed
\LZ\\
Theorem 2.4 ensures that the outcome of the Algorithm 1 is a sample from the pdf of the first-passage time $\tau_{\beta}$.\HLZ\\
To practically implement Algorithm 1, we encounter two primary challenges. The first is to find the region $B$ while the second is to define the domain $D$ for the entire trajectory of $R$.\HLZ\\
To address the first challenge, we propose choosing $B_2$ as the interval from 0 to the maximum value that $\gamma_{1}(t)+\gamma_{2}(X^{\mathbb{Q}}_t)$ can reach, which leads to the following assumption:\HLZ\\
\textbf{Assumption 2.4} $\kappa:= \sup_{t\geq0}(\gamma_1\left(t\right)+\gamma_2(X^{\mathbb{Q}}_t))$ exists.\HLZ\\

With the above assumption, we have the restricted domain $[0,\tau_{\beta}^{W}]\times\left[0,\kappa\right]$ to simulate the Poisson point process.\HLZ\\
For the second challenge, we follow \cite{herrmann2019exact} and restrict the time domain to finitely many random points. We simulate a sequence of independent random variables $\{e_k\}$; $k\geq1$, where each $e_k$ follows an exponential distribution with mean $1/\kappa$. The time variable $t$ is then defined by the cumulative sum of $e_k$ until it exceeds the time level $\tau_{\beta}^{W}$, while the space variable $u$, is defined by standard uniform random variable on the interval $\left[0, \kappa\right]$ (to compare $u$ with probability weight). Now, we can simulate $R$ on the time variable defined above. Note that this approach is one method; one can also use other alternatives mentioned in \cite{herrmann2019exact}.

For sampling $R$, we use Proposition 2.5 from \cite{herrmann2019exact}, which employs an alternative definition of the Bessel bridge (see Definition A.1 in Appendix A). This definition suggests that the Bessel bridge is the norm of a three-dimensional Brownian motion. As a result, $\eta(t)$ was formulated using a three-dimensional Brownian bridge.

Hence, we present the following new Algorithm which can be applied to practical problems:
\begin{algorithm}[h]
\caption{(Practical Implementation)}
    Let, $\beta\left(0\right)>x_{0}$.\vspace{0.1cm}\\
   \textbf{Input:} A sample of non-negative random variable $\tau_{\beta}^{W}$.\vspace{0.1cm}\\
    \textbf{Initialization}: $l=\left(0,0,0\right)$, $\mathcal{W}=0$, $\mathcal{E}_0=0$ and $\mathcal{E}_1$ an exponentially distributed random variable with mean $1/\kappa$.\vspace{0.1cm}\\
    \textbf{Step 1.} \textbf{While} $\mathcal{E}_{1}\leq \tau_{\beta}^{W}$ and $\mathcal{W}=0$,\vspace{0.1cm}\\
    \textbf{Step 2.} \hspace{0.9cm}Generate three independent random variables: a $3$-dimensional Gaussian\\ \hspace*{2cm} random variable $G$, an exponentially distributed random variable $e$ with\\\hspace*{2cm} mean $1/\kappa$ and a standard uniform random variable $u$.\vspace{0.1cm}\\
    \textbf{Step 3.} \hspace{0.93cm}Compute
        \[
        l\leftarrow\frac{\tau_{\beta}^{W}-\mathcal{E}_1}{\tau_{\beta}^{W}-\mathcal{E}_0}l+\sqrt{\frac{(\tau_{\beta}^{W}-\mathcal{E}_1)\left(\mathcal{E}_1-\mathcal{E}_0\right)}{\tau_{\beta}^{W}-\mathcal{E}_0}}G
        \]
     \hspace*{2.1cm}Update the values of $\mathcal{W}$, $\mathcal{E}_0$ and $\mathcal{E}_1$ as follows:\vspace{0.1cm}\\
    \textbf{Step 4.} \hspace{0.9cm}\textbf{if} $\kappa \cdot u \leq \gamma_1\left(\mathcal{E}_1\right)+\gamma_2\left(\beta\left(\mathcal{E}_1\right)-\|\mathcal{E}_{1}\left(\beta\left(0\right)-x_{0}\right)\left(1,0,0\right)/\tau_{\beta}^{W}+l\|\right)$ \textbf{then}\vspace{0.1cm}\\
    \textbf{Step 5.}\hspace*{1.33cm}Set $\mathcal{W}=1$\vspace{0.1cm}\\ 
    \textbf{Step 6.}\hspace*{1cm}\textbf{else}\\
    \hspace*{2.45cm}$\mathcal{W}=0$,\vspace{0.1cm}\\
    \textbf{Step 7.}\hspace*{1cm}$\mathcal{E}_0 \leftarrow \mathcal{E}_1$ and $\mathcal{E}_1\leftarrow \mathcal{E}_1 +e$.\vspace{0.1cm}\\
    \textbf{Step 8.}\hspace*{1cm}\textbf{If} $\mathcal{W}=0$ \textbf{then}\vspace{0.1cm}\\
    \textbf{Step 9.}\hspace*{1.4cm}\textbf{assign} $\tau_{\beta}=\tau_{\beta}^{W}$\vspace{0.1cm}\\
    \textbf{Step 10.} \hspace*{0.74cm}\textbf{otherwise} reject the current input and repeat the procedure with a new\\\hspace*{2.05cm} input.\vspace{0.1cm}\\
    \textbf{Step 11.} \textbf{end while}\vspace{0.1cm}\\
\textbf{Output}: A random variable $\tau$.
\end{algorithm}

\section{Numerical Assessment of the Algorithm}

In this section, we apply our Exact Algorithm and analyze its efficiency by comparing it with simulation methods that are less complex and have error rates comparable to other techniques for the type of SDE with unit diffusion described in (\ref{eq:e3}). We begin with an example in which the Exact simulation of the first-passage time of Brownian motion to a defined threshold is possible.\HLZ\\
\textbf{Example 1} Consider the following stochastic differential equation: 
\[
 dX_t=\left(K+\sin\left(X_t\right)\right)dt+dB_t,  \hspace{1cm}X_0=x_{0},
\]
and a linear threshold function $\beta(t)=at+b$, where $a$, $b$, $x_{0}$ and $K$ are real constants.\LZ\\
The SDE and the associated threshold function described above satisfy Assumption 2.2. We proceed by specifying the functions $\gamma_{1}$ and $\gamma_{2}$ relevant to this example and verifying whether these functions fulfill the necessary criteria for the application of the Exact method.
The function $\gamma_{1}:\mathbb{R}^{+}_{0}\rightarrow\mathbb{R}$ is:
\[\gamma_1\left(t\right)=-a\left(K+\sin\left(at+b\right)\right)\]
which is non-negative  $\forall t\geq0$ when $a$ and $K$ satisfy either of the two conditions: (i) $a\leq0$ and $K\geq1$, or (ii) $a\geq0$ and $K\leq-1$.\HLZ\\
Next, the function $\gamma_{2}:\mathbb{R}\rightarrow\mathbb{R}$ is: 
 \[
 \gamma_2\left(x\right)=\frac{\left({\left(K+\sin\left(x\right)\right)}^{2}+\cos\left(x\right)\right)}{2}.
 \]
Through trial and error, we observe that this function is non-negative $\forall x\in \mathbb{R}$ when $K\leq-1.6$ and $K\geq1.6$.\\
The supremum values of $\gamma_{1}$ and $\gamma_2$, denoted by $\kappa_1$ and $\kappa_2$, respectively, exist and can be expressed as follows:
\[\kappa_1=\left\{
        \begin{array}{cc}
          -a\left(K-1\right)   & ;a\geq 0 \text{ and } K\leq-1\\
          -a\left(K+1\right)   & ;a\leq 0 \text{ and } K\geq1 
        \end{array}
    \right.\]
and \[\kappa_2=\left\{
        \begin{array}{cc}
          \frac{{\left(K-1\right)}^{2}+1}{2}  & ;K\leq-1.6 \\
          \frac{{\left(K+1\right)}^{2}+1}{2}   & ;K\geq 1.6
        \end{array}
    \right. .\]
    These expressions define the supremum values under specific conditions on $a$ and $K$. Note that it is always possible to find upper and lower bounds for $\gamma_{1}$ on $[0,\tau_{\beta}^{W}]$, ensuring $\gamma_{1}$ satisfies the assumptions for all parameter values. 

Since $\gamma_{2}$ involves bounded trigonometric functions, it can also be made non-negative for $-1.6<K<1.6$  by applying a shift. Hence, for any value of the parameters, it is possible to satisfy the underlying assumptions for the method.\HLZ\\
We use the following simulation methods to analyze the efficiency of our method:
\begin{itemize}
    \item[1)] \textbf{Euler-Maruyama Method:} This method iteratively calculates the SDE path over a finite interval by discretizing the time scale and sampling the path at each time step until the threshold is reached or exceeded.\\
    Let us choose $[0,T]$ a finite time interval and discretise the time interval as $0=t_{0},t_{1},...,t_{n}=T$.\HLZ\\
    For each time step,
    \begin{list}{}{\setlength{\leftmargin}{1.5em}\setlength{\itemsep}{0.5em}}
    \item[\textit{Step 1.}] Generate the process $X$,
    \[
    X_{t_{i+1}} = X_{t_i} + \big(K + \sin(X_{t_i})\big)* (t_{i+1} - t_i) + \big(B(t_{i+1}) - B(t_i)\big)
    \]
    \item[\textit{Step 2.}] If $X_{t_{i+1}} \geq \beta(t_{i+1})$, break; and $\tau_\beta = t_{i+1}$.
\end{list}\HLZ
    The Euler-Maruyama method has $o(1)$ time-discretisation error.
    \item[2)] \textbf{Improved Euler-Maruyama Method:} The Euler-Maruyama method tends to overestimate the first-passage time. Due to time discretization, it is possible for the SDE path to cross the threshold between the time grid points without being detected by the Euler-Maruyama method. To address this issue, we implement an improvement introduced in \cite{mannella1999absorbing}, which is easy to implement. The idea is to simulate the path as in the Euler-Maruyama method and additionally check if the path has crossed the threshold between the time step $(t_{i},t_{i+1})$ with probability, \[
    P=e^{\frac{-2 d_{1} d_{2}}{(t_{i+1}-t_{i})}}
    \]
    where $d_1=d\left(X\left(t_{i}\right),\beta\left(t\right)\right),d_2=d\left(X\left(t_{i+1}\right),\beta\left(t\right)\right)$ are distance metrics. \HLZ\\
    For each time step,
    \begin{list}{}{\setlength{\leftmargin}{1.5em}\setlength{\itemsep}{0.5em}}
    \item[\textit{Step 1.}] Generate the process $X$,
    \[
    X_{t_{i+1}}=X_{t_{i}}+(K+\sin(X_{t_{i}}))*( t_{i+1}-t_{i})+(B\left(t_{i+1}\right)-B\left(t_{i}\right))
    \]
    \item[\textit{Step 2.}] $\text{Find } d_1, d_2 \text{ and } P.\text{ Also, generate a standard uniform variable } u.$
        \item[\textit{Step 3.}] $\text{If } u<P \text{ or } X_{t_{i+1}}\geq \beta\left(t_{i+1}\right),\text{ break }.$
        \item[\textit{Step 4.}] $\text{If } u<P, \text{ then } \tau_{\beta} = \left(t_{i+1}+t_{i}\right)/2,\text{ else } \tau_{\beta} = t_{i+1}$.
\end{list}\HLZ

    This method has $o(1/2)$ time-discretisation error \cite{PS_2001__5__261_0}.
\end{itemize}
To apply the Exact method, we verify that the drift function satisfies the Novikov condition. This involves checking if the following expectation is finite:
\begin{equation}
    \mathbb{E}\left[\exp\left(\frac{1}{2}\int_{0}^{\tau_{\beta}^{W}}(K+\sin(X_{s}))^{2}ds\right)\right]
    \label{eq:N1}
\end{equation}
Let us analyze the integral:
\begin{align*}
    \int_{0}^{\tau_{\beta}^{W}}(K+\sin(X_{s}))^{2}ds &= K^{2} \tau_{\beta} +  \int_{0}^{\tau_{\beta}}\sin^{2}(X_{s})ds+2K\int_{0}^{\tau_{\beta}}\sin(X_{s})ds\\
    &\leq K^{2}\tau_{\beta} + \tau_{\beta} + 2K\tau_{\beta} =(K+1)^{2}\tau_{\beta}.
\end{align*}
Since the exponential function is continuous and monotonically increasing, this directly implies:
\[
\mathbb{E}\left[\exp\left(\frac{1}{2}\int_{0}^{\tau_{\beta}^{W}}(K+\sin(X_{s}))^{2}ds\right)\right]\leq e^{\frac{(K+1)^{2}\tau_{\beta}}{2}}<\infty.
\]
We use Algorithm 2 to generate exact samples of $\tau_{\beta}$ from Example 1. We then generate approximated samples from both the Euler-Maruyama and Improved Euler methods. Finally, plot the estimated densities from the samples.
\begin{figure}[htbp] 
\centering
        \includegraphics[width=\linewidth,right]{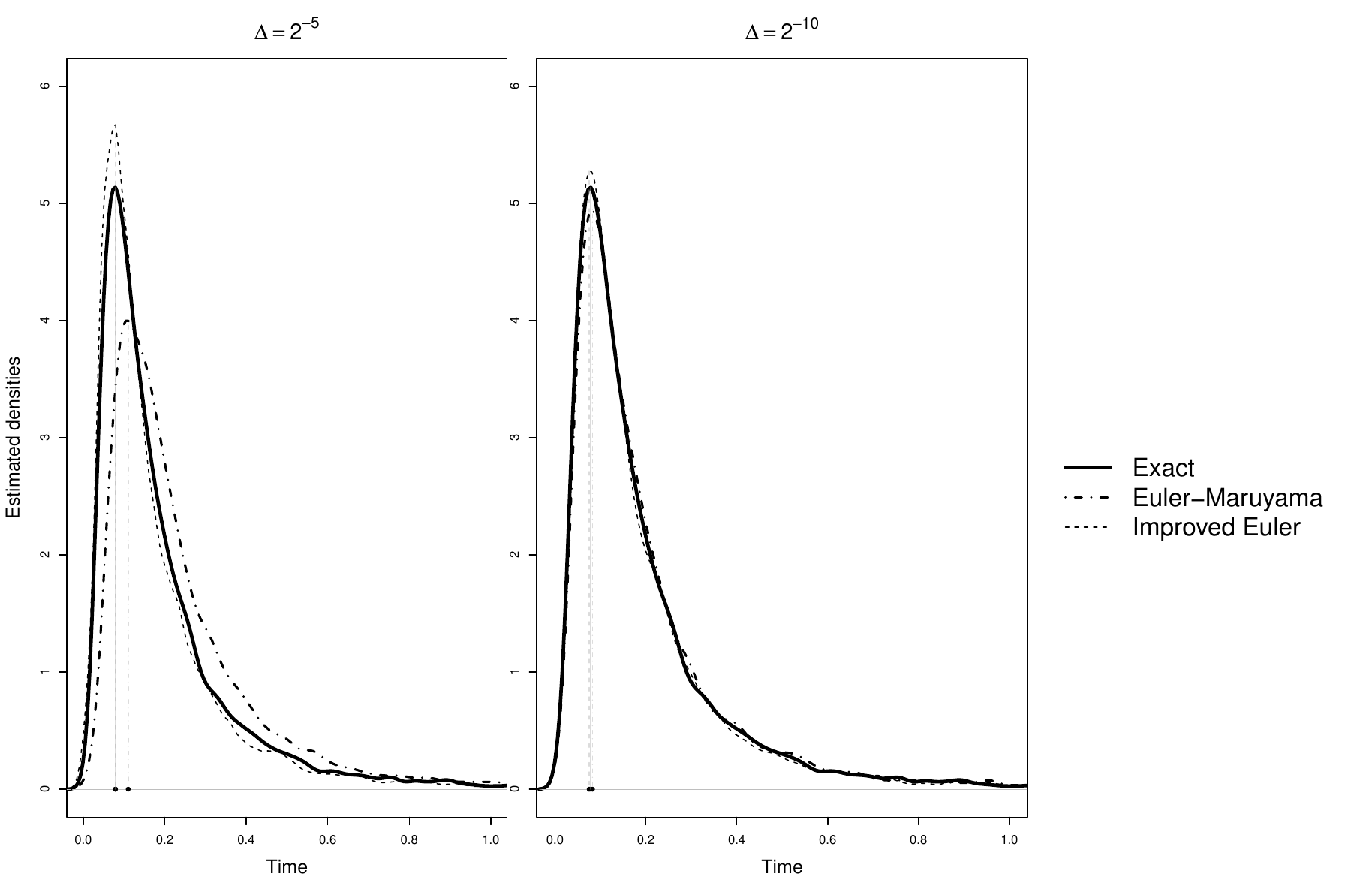}
  \caption{\footnotesize{Density plot of first-passage time for 10,000 simulations with parameters $x_{0}=0$, $K=1.6$, $a=-1$ and $b=0.5$, using three different methods, Euler-Maruyama, Improved Euler and Exact method. The points on the x-axis are the peaks of the density plots.} }
  \label{fig:Ex with linear threshold1}
\end{figure}

In Figure \ref{fig:Ex with linear threshold1}, it is evident that as the step size is reduced, the density plot from both time-discretization methods becomes closer to the density plot from the Exact method. This is expected because a smaller step size decreases the chance of time-discretization methods missing a crossing, providing a more accurate simulation. One can also observe that between the Euler-Maruyama and Improved Euler-Maruyama methods, the density from the latter is closer to the density from the Exact method. \HLZ\\
Comparing the density curves provides some insight into the accuracy of the methods. However, to ensure greater accuracy and precision, it is essential to assess "moment bias".
Moment bias refers to the discrepancy between the estimated moments of a probability distribution and the true moments (in this case, moments from the Exact method) of the underlying distribution. Specifically, we consider the first and second-moment biases, defined as:
\begin{itemize}
    \item \textit{First moment bias:} 
     \[
     \mathbb{E}\left[\hat{\tau_{\beta}}\right]-\mathbb{E}\left[\tau_{\beta}\right]
     \]
     \item \textit{Second moment bias:}
     \[
     Var\left[\hat{\tau_{\beta}}\right]-Var\left[\tau_{\beta}\right]
     \]
\end{itemize}
Here, $\hat{\tau_{\beta}}$ represents the approximated $\tau_\beta$ obtained from either the Euler-Maruyama or Improved Euler-Maruyama methods.
\begin{figure}[ht] 
  \centering     
    \includegraphics[width=\linewidth,right]{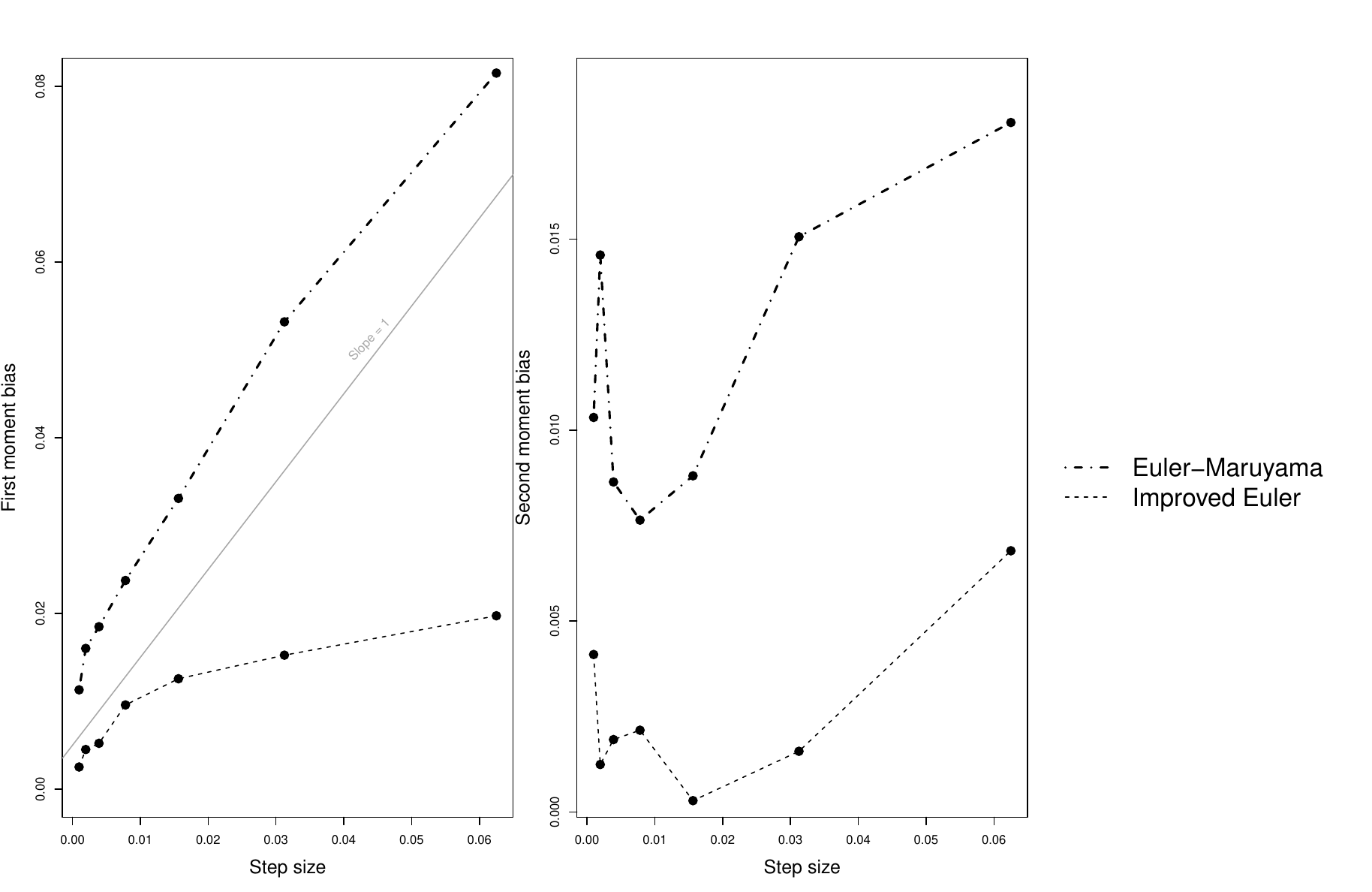}
    \caption{\footnotesize{First and second moment bias of Euler-Maruyama and Improved Euler-Maruyama method from 10k simulations with step size going from $2^{-4}$ to $2^{-10}$.} }
    \label{fig:Moment bias}
\end{figure}

In Figure \ref{fig:Moment bias}, positive bias is observed in both the first and second moments from both methods. As the step size decreases, the first moment bias decreases, but notable fluctuations persist in the second moment. This suggests the challenge of finding a balance between a smaller step size and a smaller statistical error when using a time-discretization method. One can notice that even in terms of statistical precision, the Improved Euler method performs better than the Euler method.\HLZ\\
Although we know that the Exact method gives accurate results, and if we want to achieve better accuracy in time-discretization methods, it is difficult to maintain statistical consistency. However, to say that the Exact method is efficient to use depends on other factors as well, such as how much time it takes to run. If the Exact method takes significantly longer, we must consider whether the gain in accuracy justifies the increased computational expense.\HLZ\\
To test efficiency accordingly, we compare the computational time of the Exact method with that of the Euler-Maruyama and Improved Euler-Maruyama methods and conduct the Kolmogorov-Smirnov test. In this way, we can simultaneously compare how much time the extra accuracy costs. 
This analysis will provide a more complete understanding of the trade-offs between accuracy and computational efficiency.
\begin{table}[ht]
\label{tab:1}       
\caption{KS-test results and computational time from $10^4$ simulations.}
\centering
\begin{tabular}{llllll}
\hline\noalign{\smallskip}
$\Delta$&\multicolumn{2}{l}{Euler-Maruyama} & \multicolumn{2}{l}{Imp. Euler-Maruyama} & Exact  \\
\noalign{\smallskip}\hline\noalign{\smallskip}
 & p-value &Time(s) & p-value &Time(s)& Time(s) \\
\noalign{\smallskip}\hline\noalign{\smallskip}
$2^{-4}$&  $0$ & $0.53$  & $0$ & $0.67$ &  \multirow{7}{*}{$11.97$}\\
         $2^{-5}$& $0$ & $0.66$ & $0$ & $0.9$ &  \\
         $2^{-6}$& $0$ & $0.93$ & $0$ & $1.19$ & \\
         $2^{-7}$& $0$ & $1.34$ & $1.264\times10^{-12}$ & $2$ & \\
         $2^{-8}$& $0$ & $2.34$ & $1.665\times10^{-5}$ & $3.51$ & \\
         $2^{-9}$& $2.859\times10^{-7}$ & $4.56$ &  $0.0002324$ & $6.28$ & \\
        $2^{-10}$& $0.0008803$ & $8.35$ &  $0.02755$ & $12.14$ & \\
\noalign{\smallskip}\hline
\end{tabular}
\caption*{\footnotesize \textbf{Note}: The Exact method does not involve time-discretisation error. Thus, in any case of $\Delta$, we have the same computational time from the Exact method, which here is 11.97 seconds.}
\end{table}

Table 1 shows us when the time-discretization methods start getting closer to a good p-value. We see that the density from Euler-Maruyama, even with $2^{-10}$ step size, gives a very small p-value. The Improved Euler-Maruyama method with the same step size gives a p-value close to 0.05 but with more computational time than the Exact method. Thus, we observe that both time-discretization methods require more computational time than the Exact method to achieve comparable accuracy.\HLZ\\
   It is important to note that the time complexity of the methods is influenced by various parameters. The performance of the Exact Algorithm is notably sensitive to $\kappa$, $a$ and $b$. For instance, setting $K>1.6$ increases the value of $\kappa$, which directly affects the algorithm's processing time. For example, with $K = 2$, the computational time for the Exact method from 10,000 simulations becomes 16 seconds. Additionally, higher values of $a$ and $b$ push the threshold away, resulting in an increase in $\tau_{\beta}^{W}$ and thereby the while loop in the algorithm will run for a longer duration. This will be further discussed in Section 4.
   
\subsection{Estimating the unknown Distribution of the FPT of Brownian Motion }
In this section, we discuss the case of threshold functions for which the distribution of the first-passage time of Brownian motion is not available. In such cases, employing an approximation technique with robust convergence properties is beneficial. We propose using the iterative approach introduced by Hermann and Tanr{\'e} \cite{herrmann2016first}. This method is particularly interesting due to its strong convergence properties, allowing for adjustments in algorithm parameters to enhance accuracy without significantly impacting computational time.\HLZ\\
Let us discuss the process of obtaining the first-passage time of Brownian motion using this iterative approach. In scenarios where the threshold function, $\beta(t)$, is non-decreasing and $\beta\left(0\right)>0$, the intuition is that the Brownian motion needs to successively cross a sequence of horizontal lines before reaching the threshold. We commence with $\beta\left(0\right)$ and determine the first-passage time of Brownian motion to $\beta(0)$, denoted by $\tau_{1}$. Next, we determine the first-passage time to  $\beta(\tau_{1})$, denoted by $\tau_{2}$, and continue this process iteratively. We stop the iteration when the difference $\beta(\tau_{n})
-\beta(\tau_{n-1})$ becomes very small (less than a very small parameter $\epsilon$). At this point, $\tau_{n}$ is considered to be the approximated first-passage time for the given $\epsilon$. Let $\{\tau_{\epsilon}^{W}\}_{\epsilon}$ be a sequence of approximated first-passage time with different $\epsilon$. As $\epsilon$ tends to 0, this sequence converges to $\tau_{\beta}^{W}$.

In cases where the threshold function is curvy, i.e., a time-dependent function that is not non-decreasing,  we can effectively resolve the situation by tilting successive imaginary horizontal lines. This tilting is performed so that the common slope corresponds to $r\leq\inf_{t\geq 0} \beta' \left(t\right)$. Subsequently, we obtain the sequence by performing similar iterations, but this time using successive linear straight lines.\HLZ\\
Thus, for $\epsilon$ small enough, $\tau_{\epsilon}^{W}$ is a good approximation of $\tau_{\beta}^{W}$, we can then apply the Exact Algorithm using the rejection probability as follows: \[
\eta^{\epsilon}(t)=\mathbb{E}_{\mathbb{Q}}\left[ \exp\left( -\int_{0}^{\tau_{\epsilon}^{W}}(\gamma_{1}(s) + \gamma_{2}(X_s^{\mathbb{Q}}))ds \right) |t=\tau_{\epsilon}^{W}\right].
\]
Since $\gamma_1\left(t\right)+ \gamma_2(X^{\mathbb{Q}}_t)$ is continuous on $[0,\tau_{\beta}^{W}]$, and $\tau_{\beta}^{W}$ is positive, we can conclude that
\[
\eta^{\epsilon}(t) \longrightarrow \eta(t).
\]
Then, the Girsanov transformation from equation (\ref{eq:e5}) becomes 
\begin{align*}
    \mathbb{E}_{\mathbb{P}}[\Psi(\tau_{\beta}) \mathbf{1}_{\{\tau_{\beta}<\infty\}}] = \mathbb{E}_{\mathbb{Q}}[\Psi(\tau_{\beta})\eta(t)] = \mathbb{E}_{\mathbb{Q}}[\Psi(\tau_{\beta})\lim_{\epsilon\to 0}\eta^{\epsilon}(t)] =\lim_{\epsilon \to 0}\mathbb{E}_{\mathbb{Q}}[\Psi(\tau_{\beta})\eta^{\epsilon}(t)].
\end{align*}
Thereby, the last equality holds due to the Dominated Convergence theorem, as $\eta^{\epsilon}$ is convergent and bounded.\HLZ\\
Let us now discuss an example with a threshold function for which the first-passage time of Brownian motion cannot be exactly simulated.\LZ\\
\textbf{Example 2} 
Consider again the SDE 
\[
 dX_t=\left(K+\sin\left(X_t\right)\right)dt+dB_t,   \hspace{1cm}X_0=x_{0},
\]
but now with an exponential threshold $\beta(t) = ae^{-bt}$.\LZ\\
Then
\[\gamma_1\left(t\right)=abe^{-bt}\left(K+\sin\left(ae^{-bt}\right)\right).\]
The function $\gamma_1$ is non-negative $\forall t \geq 0$ when $a,b \text{ and } K$ satisfies either of the two conditions, (i) $a\cdot b\geq 0$ and $K\geq 1$, or (ii) $a\cdot b\leq 0$ and $K\leq -1$.
Since $\gamma_1$ is a decreasing function when $b\geq 0$, its supremum must occur at $t=0$. When $b\leq 0$, $\gamma_{1}$ is increasing, and thus, its supremum does not exist on $\mathbb{R}$. Note that we only need $\gamma_{1}$ to satisfy the conditions on $[0,\tau_{\beta}^{W}]$, which is always possible as $\gamma_{1}$ is continuous and also deterministic for a particular $\tau_{\beta}^{W}$. \HLZ\\
Now consider,
\[
\gamma_2\left(x\right)=\frac{\left({\left(K+\sin\left(x\right)\right)}^{2}+\cos\left(x\right)\right)}{2}.
\]
We obtain non-negativity and its supremum as in Example 1.\HLZ\\
The threshold function, $\beta\left(t\right)$ is, however, not increasing, so we use tilted horizontal lines to approximate the $\tau_{\beta}^{W}$, cf. \cite{herrmann2016first}.\HLZ\\
We generate $\tau_{\beta}^{W}$ using the iterative algorithm described in Appendix B., and then apply Algorithm $2$ to accept the samples for the first-passage time of $X_t$ to $\beta(t)$, resulting in the following outcomes:

\begin{figure}[ht]
\centering
    \includegraphics[scale=0.43]{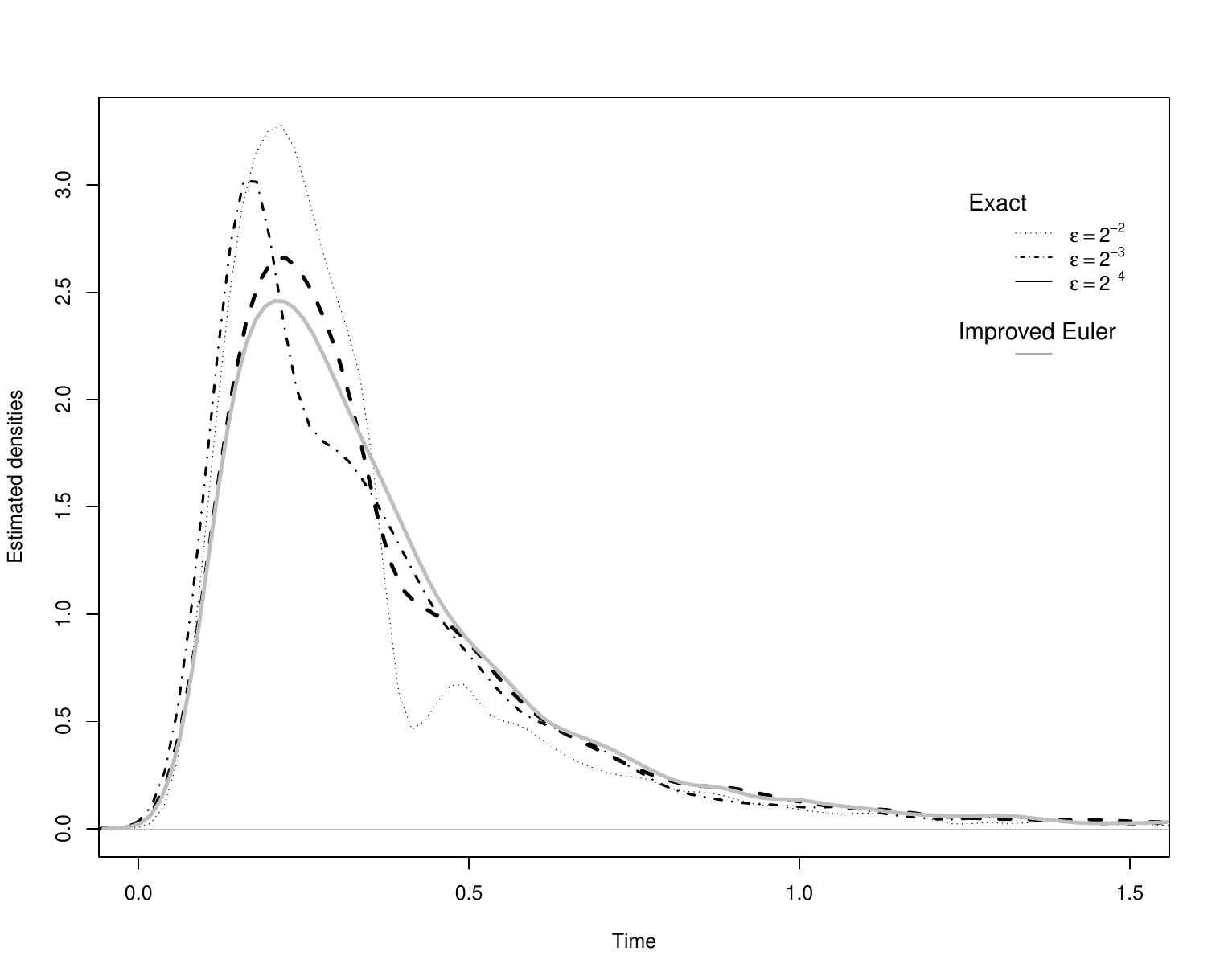} 
    \caption{\footnotesize{Density plot of first-passage time from 5000 simulations with parameters $x_{0}=0$, $K=1.6$, $a=b=1$ using Exact method (with different $\epsilon$ parameter for generating $\tau_\beta$) and "High-resolution" Improved Euler method.}}
    \label{fig:epsilon density plot}
\end{figure}
Figure \ref{fig:epsilon density plot} clearly shows that with smaller $\epsilon$, density curves using the Exact method get closer to the density curve using the High-resolution Improved Euler method, which is a good sign. However, relying solely on this observation is insufficient to ascertain the method's quality. So, let us analyze it further by checking the bias and accuracy using statistical tests and the computational time of the algorithm.
\begin{table}[ht]
\label{tab:2}       
\caption{KS-test results and computational time from $5\cdot10^3$ simulations.}
\centering
\begin{tabular}{lllll}
\hline\noalign{\smallskip}
& \multicolumn{3}{l}{Exact} & High-resolution\\
& & & &Improved Euler \\
\noalign{\smallskip}\hline\noalign{\smallskip}
$\epsilon$ & $2^{-2}$ & $2^{-3}$ & $2^{-4}$ & \\
\noalign{\smallskip}\hline\noalign{\smallskip}
 Mean & 0.36828 & 0.40583 & 0.45055 & 0.43931 \\
    Variance & 0.20811 & 0.23825 & 0.33214 & 0.23156 \\
    Time to run (s) & 27.17 & 30.41 & 33.68 & \\
    p value & 0 & 0  & $0.4247$ &\\
\noalign{\smallskip}\hline
\end{tabular}
\end{table}\\
We deduce from Table 2 that with smaller values of $\epsilon$, the accuracy notably improves, as indicated by a lower bias and higher p-values. Concerning time complexity, a decreasing $\epsilon$ results in an increase in the run-time. However, the extent of this increase is moderate. This is primarily because changes in $\epsilon$ predominantly impact the time required for generating $\tau_{\beta}^{W}$, while the parameters influencing the run-time of Algorithm 2 remain constant. Additionally, prior research detailed in \cite{herrmann2016first} also emphasizes the marginal impact of $\epsilon$ on the time needed for generating $\tau_{\beta}^{W}$.\LZ\\
\begin{remark}
    \begin{itemize}
    \item[1.] Even in this Example, $\kappa$ and $\tau_{\beta}^{W}$ directly correlate with the run-time. Therefore, an increase in $K$ beyond $1.6$ or any other factor causing an escalation in $\kappa$, will increase the run-time. \\
    Similarly, a larger $a$ or a smaller $b$, or other factors leading to a rise in $\tau_{\beta}^{W}$, will similarly contribute to an increase in run-time.
    \item[2.] The parameter $r$ affects the run-time of $\tau_{\beta}^{W}$ as demonstrated in Theorem 2.1 of \cite{herrmann2016first}. Different values of $r$ yield different results, but there is no consistent pattern. Therefore, it might be feasible to choose the smallest possible $r$ and adjust $\epsilon$ to achieve good accuracy.
\end{itemize}
\end{remark}
\section{The Time Complexity of the Algorithm}

We discussed the Exact Algorithm in the previous sections and compared its performance with other simulation methods. While alternative methods can attain comparable accuracy by minimizing step sizes, they often introduce a bias. However, it is crucial to note that the Exact Algorithm may require more computational time compared to other methods. In this section, we will estimate the upper bound of the time complexity of the Exact Algorithm. If this bound is found to be high, we will investigate potential strategies to optimize and reduce it.\HLZ\\
Let $I$ be the number of iterations observed to simulate $\tau_\beta$ using the Exact Algorithm described in Section 2.\\
The expected number of iterations until we accept an outcome $\tau_{\beta}$ depends on the probability of acceptance. This relationship, combined with Theorem 2.2, yields the following equality:
\[\mathbb{E}\left[I\right]=\mathbb{E}_{\mathbb{Q}}\left[\eta\left(\tau_{\beta}\right)\right]^{-1}=\exp\{-A(x_{0})+A\left(\beta(0)\right)\}.\]
The function $\gamma_1(t)+\gamma_2(X^{\mathbb{Q}}_{t})$ is bounded by $\kappa$, leading to the inequality\\
\begin{equation}
    \mathbb{E}_{\mathbb{Q}}\left[\eta\left(\tau_{\beta}\right)\right] \geq \mathbb{E}_{\mathbb{Q}}\left[e^{-\kappa \cdot \tau_{\beta}}\right].
    \label{eq:eq8}
\end{equation}
Clearly, the number of iterations depends on $\kappa$ and the first-passage time of Brownian motion. We also observed the same effect in the examples in Section 3. The bound for the number of iterations is exponential, meaning that larger values of $\kappa$ and $\tau_{\beta}^{W}$ can increase the algorithm's run time exponentially. Following this discussion, we will introduce some approaches to reduce this run time.

\subsection{Shifting $\gamma_1(t)+\gamma_2(X^{\mathbb{Q}}_{t})$}
Let us assume that there exists a lower bound of $\gamma_1(t)+\gamma_2(X^{\mathbb{Q}}_t)$ on $\left[0,\infty\right)$.\\
The infimum of the sum of two non-negative functions should be equal to the sum of the infimum of the individual functions.\HLZ\\
Let, $\gamma_{0,1}=\underset{t\geq0}{\inf} \gamma_1\left(t\right)$ and $\gamma_{0,2}=\underset{x\in \mathbb{R}}{\inf} \gamma_2\left(x\right)$, then\\
\[\gamma_{0,1}+\gamma_{0,2}\leq\underset{t\geq0}{\inf}\left(\gamma_1\left(t\right)+\gamma_2(X^{\mathbb{Q}}_t)\right)\leq \kappa_1 +\kappa_2.\]
Thus, the rejection probability can be computed using the shifted function $\gamma_1(t)+\gamma_2(W_{t})-\gamma_{0,1}-\gamma_{0,2}$. Consequently, the corresponding upper bound becomes the shifted value $\kappa_1+\kappa_2-\gamma_{0,1}-\gamma_{0,2}$, leading to a reduced number of iterations.

If only one of the functions has a lower bound, shifting that function alone will reduce $\kappa$ and consequently the run time as well.
\subsection{Space Splitting}
The idea is to split the space between the initial value and the threshold function by constructing $k$ functions, each obtained by shifting the threshold function. Let $\beta_i\left(t\right);1\leq i\leq k$ be the sequence of these functions.
We then find the first-passage time, $\tau_{\beta_{i}}$ of $X\left(t\right)$ to each of the functions, $\beta_i(t)$.\HLZ\\ The rejection probability and the Exact Algorithm to generate $\tau_{\beta_{1}}$ are the same as described in Section 2. For $\tau_{\beta_{i}};2\leq i \leq k$, we use $\tau_{\beta_{i-1}}$ (already accepted), as described in the following subsection.

\subsubsection{Rejection Probability} 
We know that the Girsanov transformation for the first-passage time of $X\left(t\right)$ to $\beta_1\left(t\right)$ using Theorem 2.2 is given as follows: 
\begin{equation}
\begin{split}
    \mathbb{E}_{\mathbb{P}}[\Psi \left(\tau_{\beta_1}\right)&\mathbf{1}_{\{\tau_{\beta_1}<\infty\}}]= \mathbb{E}_\mathbb{Q}\left[\Psi \left(\tau_{\beta_1}\right)e^{-A\left(x\right)}\exp\left(A(\beta_{1}({\tau}_{\beta_1}))-\int_{0}^{{\tau}_{\beta_1}} {\gamma}_{2}(X^{\mathbb{Q}}_s)ds\right)\right].
\end{split}
\label{eq:eq9}
\end{equation}

\noindent We now discuss the rejection probability to obtain the first-passage time of the same path of $X(t)$ to $\beta_2\left(t\right)$. Using the Girsanov transformation for $\tau_{\beta_{2}}$, we get \LZ\\

\begin{align*}
\mathbb{E}_{\mathbb{P}}[\Psi (\tau_{\beta_2})\mathbf{1}_{\{\tau_{\beta_2}<\infty\}}]&= \mathbb{E}_{\mathbb{Q}} \left[\Psi({\tau}_{\beta_{2}}) \exp\left( \int_{0}^{\tau_{\beta_{2}}} \alpha(X^{\mathbb{Q}}_{s})dB_{s} -\frac{1}{2} \int_{0}^{\tau_{\beta_{2}}} \alpha^{2}(X^{\mathbb{Q}}_{s})ds \right)\right]\\
&= \mathbb{E}_{\mathbb{Q}}\left[\Psi(\tau_{\beta_2})e^{-A\left(x\right)} \exp\left(A(\beta_{2}({\tau_{\beta_{2}}})) -\int_{0}^{\tau_{\beta_{2}}} \gamma_{2}(X^{\mathbb{Q}}_s)ds\right)\right]\\
&= \mathbb{E}_{\mathbb{Q}}\left[\Psi(\tau_{\beta_2})e^{-A(x)} \exp\Biggl(A(\beta_{1}(\tau_{\beta_1}))-A(\beta_{1}(\tau_{\beta_1}))+A(\beta_{2}(\tau_{\beta_{1}}))\right.\Biggr.\\ 
&\quad \left.\left.-A(\beta_{2}(\tau_{\beta_{1}}))+A(\beta_{2}(\tau_{\beta_{2}}))-\int_{0}^{\tau_{\beta_{1}}} \gamma_{2}(X^{\mathbb{Q}}_s)ds-\int_{\tau_{\beta_{1}}}^{\tau_{\beta_{2}}} \gamma_{2}(X^{\mathbb{Q}}_s)ds\right)\right]\\
    &=\mathbb{E}_{\mathbb{Q}}\Biggl[\Psi(\tau_{\beta_{2}})\exp\Biggl(-A(\beta_{1}(\tau_{\beta_{1}}))+A(\beta_{2}(\tau_{\beta_{1}}))-  \int_{\tau_{\beta_{1}}}^{\tau_{\beta_{2}}} \gamma_{1}(X^{\mathbb{Q}}_s) ds\\ 
    &\quad-\int_{\tau_{\beta_{1}}}^{\tau_{\beta_{2}}} \gamma_{2}(X^{\mathbb{Q}}_s)ds\Biggr)\Biggr] \times \frac{{\mathbb{E}_{\mathbb{Q}}\left[\Psi\left(\tau_{\beta_{1}}\right)\right]}}{\mathbb{E}_{\mathbb{P}}\left[\Psi(\tau_{\beta_{1}})\right]},
\end{align*}
where in the last line (\ref{eq:eq9}) was used.

Since, $\tau_{\beta_{1}}$ is the one that we already accepted for threshold function, $\beta_{1}(t)$, which means $\tau_{\beta_{1}}$ is known at this stage.\HLZ\\
It also implies that ${\mathbb{E}_{\mathbb{Q}}\left[\Psi\left(\tau_{\beta_{1}}\right)\right]}=\mathbb{E}_{\mathbb{P}}\left[\Psi(\tau_{\beta_{1}})\right]$ . And hence the equality becomes\HLZ\\
$\mathbb{E}_{\mathbb{P}}[\Psi (\tau_{\beta_2})\mathbf{1}_{\{\tau_{\beta_2}<\infty\}}]$
\[=\mathbb{E}_{\mathbb{Q}}\left[\Psi\left(\tau_{\beta_{2}}\right)e^{-A(\beta_{1}(\tau_{\beta_{1}}))+A(\beta_{2}(\tau_{\beta_{1}}))}\exp\left(- \int_{\tau_{\beta_{1}}}^{\tau_{\beta_{2}}} \gamma_{1}(X^{\mathbb{Q}}_s) -\int_{\tau_{\beta_{1}}}^{\tau_{\beta_{2}}} \gamma_{2}(X^{\mathbb{Q}}_s)ds\right)\right]\]

The algorithm to generate $\tau_{\beta_{2}}$ using the above rejection probability would be the same as Algorithm 2 with an updated Bessel process, \[
R\left(s\right)=\beta_{2}\left(\tau_{\beta_{1}}+\tau_{\beta_{2}}-s\right)-X^{\mathbb{Q}}_{\tau_{\beta_{1}}+\tau_{\beta_{2}}-s}, \hspace{0.5cm}\tau_{\beta_{1}}\leq s \leq \tau_{\beta_{2}},
\] starting at 0 and ending at $\beta_{2}\left(\tau_{\beta_{1}}\right)-\beta_{1}\left(\tau_{\beta_{1}}\right)$.
\subsubsection{Efficiency of the Algorithm using Space Splitting}
Let us assume that the threshold function  
$\beta\left(t\right)$ is linear with $a\cdot b <0$. In this case, $\tau_{\beta}^{W}$ follows an Inverse Gaussian distribution for which the moment-generating function exists. Therefore, inequality (\ref{eq:eq8}) becomes:
\begin{align*}
\mathbb{E}_{\mathbb{Q}}\left[\eta\left(\tau_{\beta}\right)\right] \geq \mathbb{E}\left[e^{-\kappa \cdot \tau_{\beta}}\right] =\exp{\left(-ab+ab\sqrt{1+\frac{2\kappa}{a^2}}\right)}\,,
\end{align*}
which implies
\[\mathbb{E}_{\mathbb{Q}}\left[I\right]\leq\exp{\left(ab-ab\sqrt{1+\frac{2\kappa}{a^2}}\right)}.\]
Note that the bound on $\mathbb{E}_{\mathbb{Q}}\left[I\right]$ is exponential.\HLZ\\
Now let $I_{\text{shift}}$ be the number of iterations observed when using Space Splitting. Then the upper bound of $\mathbb{E}_{\mathbb{Q}}\left[I_{\text{shift}}\right]$ is the sum of upper bounds of expected iterations to generate $\tau_{\beta_{i}}$ for each $i$.
\begin{align*}
    \mathbb{E}\left[I_{shift}\right]&\leq \sum_{i=1}^{k} \exp{\left(a\cdot{b}_{i} \left(1-\sqrt{1+\frac{2\kappa_i}{a^2}}\right)\right)}\\
    &\leq \sum_{i=1}^{k}\exp{\left(a\cdot\frac{b}{k}\left(1-\sqrt{1+\frac{2\max_{i}{\kappa_i}}{a^2}}\right)\right)}\\
    &=k\exp{\left(\frac{ab}{k}\left(1-\sqrt{1+\frac{2\max_{i}{\kappa_i}}{a^2}}\right)\right)}\\
\end{align*}
By choosing $k=\lfloor ab\left(1-\sqrt{1+\frac{2\max_{i}{\kappa_i}}{a^2}}\right)\rfloor +1$, the inequality becomes
\[
\mathbb{E}\left[I_{shift}\right]
\leq \left(\lfloor ab\left(1-\sqrt{1+\frac{2\max_{i}{\kappa_i}}{a^2}}\right)\rfloor +1\right)e.
\]
The bound now becomes linear. Similarly, we can show bound to 
be linear for the case when $a\cdot b >0$.\\
A similar result was proved in \cite{herrmann2016first} (Section 3.1.1) for the case when the threshold $\beta(t)$ is a constant.\HLZ\\
If we use the algorithmic approach from \cite{herrmann2016first} as described in Section 3.1 to generate the first-passage time of Brownian motion, we use the first-passage time to either constant or linear functions to approximate the first-passage time to non-decreasing and other curvy functions, respectively.\\
Assume that we have a non-decreasing threshold and we approximate $\tau_{n}=\tau_{\epsilon}^{W}$ as $\tau_{\beta}^{W}$. This means that $\tau_{n}$ is the FPT of Brownian motion to a constant $\beta(\tau_{n-1})$ such that $|\beta(\tau_{n-1})-\beta(\tau_{n-2})|<\epsilon$. Using the same approximation technique, we approximate the first-passage time to $\beta_{i}(t)$, $\tau_{n}^{[i]}={\tau_{\epsilon}^{W}}^{[i]}$. This means that $\tau_{n}^{[i]}$ is the FPT of Brownian motion to a constant $\beta_{i}(\tau_{n-1}^{[i]})$ such that $|\beta(\tau_{n-1})-\beta(\tau_{n-2})|<\epsilon$. Instead of thinking space splitting of threshold $\beta(t)$, we can think of it as space splitting of constant threshold $\beta(\tau_{n})$ with k constant thresholds $\beta_{i}(\tau_{n}^{[i]}); 1\leq i\leq k$. We moreover know that for constant thresholds, with suitable $k$, the bound becomes linear. Similarly, when we have a curvy threshold, instead of thinking of space splitting of threshold $\beta(t)$, we can think of it as space splitting of the tilted horizontal line. 
\section{Predicting Spike Times in a Neuron using Exact simulation}
\begin{figure}[ht]
        \centering
        \includegraphics[width=0.8\linewidth]{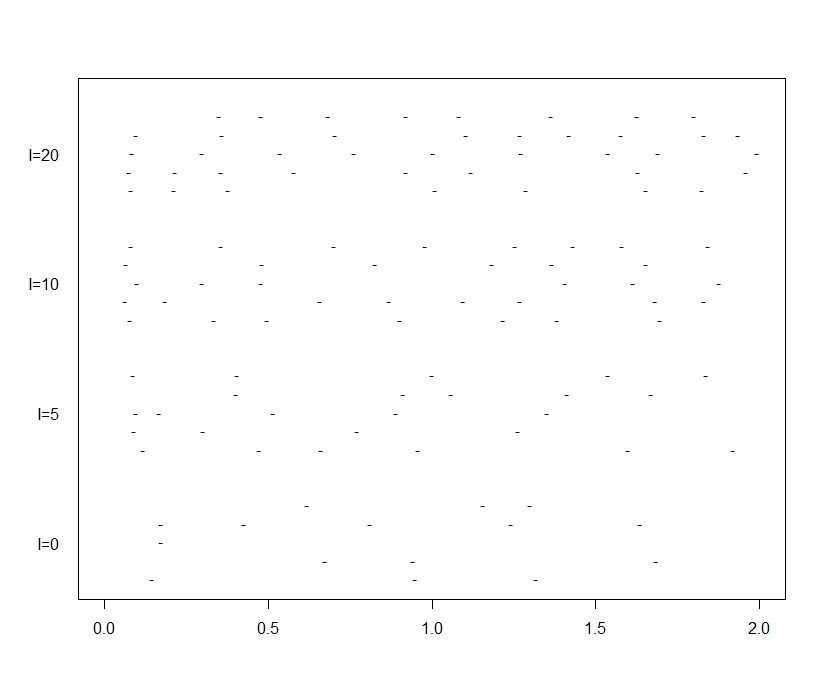}
        \caption{\footnotesize Spike times of a neuron with different input current, $I$, observed until 2 seconds. For each $I$, we did 5 trials with parameters $\tau=-1$, $V_r =1$, $\sigma=1$, $v_0 =1$, $\theta_0 =1$, $\Delta=1$ and $\tau_1 =1$.}
        \label{fig:enter-label}
\end{figure}

We conclude our discussion with an application to a neuron model. A neuron generates a spike when its membrane voltage reaches a specific threshold, and a spike train is the sequence of recorded times at which these spikes occur. The pattern of this spike train conveys information to the nervous system.\HLZ\\
Standard models, such as the linear "leaky integrate-and-fire" model with a constant threshold, are often insufficient to capture the firing properties of real neurons (cf. \cite{gerstner2014neuronal}). To address this, we consider a more complex, non-linear leaky integrate-and-fire model to describe the membrane voltage of a neuron, represented using the following SDE:
\begin{equation}
    dV_t=(f(V_t) + R(V_t)\cdot I)dt - \sigma V_tdB_{t}; \hspace{0.4cm} V_0=v_{0}, \hspace{0.2cm} \sigma>0
    \label{eq:eq10}
\end{equation}
where $I$ is the input current, $f$ is the leak term and $R$ is input resistence.\HLZ\\
In \cite{levakova2019adaptive}, the authors observed that using an adaptive exponentially declining threshold for detecting spike times is highly effective, as it allows the spike times to align more closely with real data. We describe this adaptive threshold as follows:

\begin{itemize}
    \item[1.] When the neuron does not generate a spike, $\theta(t)$ exponentially declines to $\theta_{0}>v_{0}$, i.e.
    \begin{equation}
        \tau_{1} \frac{d\theta(t)}{dt} = -(\theta(t)-\theta_{0}).
        \label{eq:eq11}
    \end{equation}
    \item[2.] To add the adaption variable, we adjust $\theta(t)$ as follows:
    \begin{equation}
        \theta(t) = \theta_{0} + (\theta(\tau^{+}) - \theta_{0})\exp\left(-\frac{t-\tau}{\tau_{1}}\right), \hspace{0.4cm} \forall t\geq \tau,
        \label{eq:eq12}
    \end{equation}
    where $\tau$ is the time of the last spike and $^{+}$ represents the limit from above.\\
    If the voltage reaches the threshold at time $\tau$, $V_{\tau}\geq\theta(\tau)$, the threshold increases by a step $\Delta/\tau_{1}$, i.e.
    \begin{equation}
        \theta(\tau^{+}) = \theta(\tau^{-}) + \Delta/\tau_{1},
        \label{eq:eq13}
    \end{equation}
    where $\Delta$ represents the strength of adaption due to a single spike.
\end{itemize}
How does this model work? The membrane voltage $V$ evolves according to the SDE described in equation (\ref{eq:eq10}). The evolution of the membrane voltage continues until it reaches a threshold  \[\theta_{1}(t) = \theta_{0} + \exp(-t/\tau_{1}).\] Say $t_1$ is the time at which the voltage reaches the threshold $\theta_{1}$. At this point, a spike occurs, and the membrane voltage resets to its resting position. After the reset, the membrane voltage evolves again according to the same SDE, and the process repeats when the voltage reaches the new shifted threshold, i.e. 
\[\theta_{2}(t) = \theta_{0} + (\theta_{1}(t_{1}) - \theta_{0} + \Delta/\tau_{1})\exp(-(t-t_{1})/\tau_{1}).
\] Say $t_2$ is the next spike time. Again, the membrane voltage resets to its resting position after this spike, and the process continues.\HLZ\\
Thus, a trial of spike train will be the sequence of spike times $\{t_1,t_2,...\}$. \HLZ\\
We now consider the quadratic leaky integrate-and-fire model, such that 
\[
f(x) = -\frac{1}{\tau}x(x-V_{r}),\hspace{0.2cm}\text{and}\hspace{0.2cm}R(x) = x,
\]
where $V_{r}$ is the resting position and $\tau$ is integration time constant.\HLZ\\
We set $X(t) = -\frac{1}{\sigma}\ln(V(t))$ to transform the SDE to an SDE with unit diffusion. Using the transformation, we obtain SDE for $X(t)$ as below:
\begin{equation}
    dX(t) = \left( \left( \frac{\sigma}{2} - \frac{I}{\sigma} - \frac{V_r}{\tau \sigma} \right) + \frac{1}{\tau \sigma}e^{-\sigma x} \right)dt + dB_{t}, \hspace{0.4cm} X(0) = -\frac{1}{\sigma}\ln(V_{0}).
    \label{eq:eq14}
\end{equation}
As mentioned in Section 5 of \cite{beskos2006retrospective}, the process $V_t$ won't hit 0 with probability 1, and thus the transformation is well-defined.\HLZ\\
The time at which the SDE described in equation (\ref{eq:eq10}) hits the threshold $\theta(t)$, the same time SDE in equation (\ref{eq:eq14}) hits $\beta(t) = -\frac{1}{\sigma}\ln(\theta(t))$.\HLZ\\
The associated $\gamma_{1}$ and $\gamma_{2}$ functions are then:
\begin{equation}
    \gamma_{1}(t) = -\left( \left( \frac{\sigma}{2} - \frac{I}{\sigma} - \frac{V_{r}}{\tau \sigma} \right) + \frac{1}{\tau \sigma}e^{-\sigma \beta(t)} \right)\times \beta'(t)
    \label{eq:eq15}
\end{equation}
\begin{equation}
    \gamma_{2}(t) = \frac{1}{2}\left[ \left( \frac{\sigma}{2} - \frac{I}{\sigma} - \frac{V_r}{\tau \sigma} \right)^{2} + \frac{1}{\tau^{2} \sigma^{2}}e^{-2\sigma x} - \frac{2}{\tau \sigma}\left( \frac{I}{\sigma} + \frac{V_r}{\tau \sigma}e^{-\sigma x} \right) \right]
    \label{eq:eq16}
\end{equation}
The SDE in equation (\ref{eq:eq14}) and associated threshold $\beta(t)$ satisfies the necessary assumptions to use the Exact Algorithm.\HLZ\\
However, one can observe that the $\gamma_2$ function is unbounded at $-\infty$. We use the fact that we apply the Exact Algorithm on the interval $\left[0,\tau_{\beta}\right]$. In this interval, the Brownian motion can go minimum to $\beta(t)$, since $X(0)>\beta(0)$. This implies that on the interval $\left[0,\tau_{\beta}\right]$, $\gamma_{2}$ is bounded.\\
To apply the Exact method, we verify whether the drift function satisfies the Novikov condition:
\[
\mathbb{E}\left[\exp\left(\frac{1}{2}\int_{0}^{\tau_{\beta}}\left(\frac{\sigma}{2}-\frac{I}{\sigma}-\frac{V_{r}}{\tau \sigma}+\frac{1}{\tau \sigma}e^{-\sigma X_{s}}\right)^{2}\right)\right]<\infty.
\]
Let us analyze the integral:\HLZ\\
\[\int_{0}^{\tau_{\beta}}\left(\frac{\sigma}{2}-\frac{I}{\sigma}-\frac{V_{r}}{\tau \sigma}+\frac{1}{\tau \sigma}e^{-\sigma X_{s}}\right)^{2}ds.\]
This can be expressed as follows:
\begin{align*}
    \left(\frac{\sigma}{2}-\frac{I}{\sigma}-\frac{V_{r}}{\tau \sigma}\right)^{2}\tau_{\beta}+\frac{1}{\tau^{2}\sigma^{2}}\int_{0}^{\tau_{\beta}}e^{-2\sigma X_{s}}ds+\frac{2}{\tau \sigma}\int_{0}^{\tau_{\beta}}e^{-\sigma X_{s}}ds
\end{align*}
Since the threshold $\theta(t)$ is above the process $V_{t}$, we have:
\begin{align*}
    V_{t}&\leq \theta(t) \text{ on } [0,\tau_{\beta}] \hspace{0.4cm}
    \Rightarrow -\frac{1}{\sigma}\ln(V_{t})\geq -\frac{1}{\sigma}\ln(\theta(t)) \hspace{0.4cm}
    \Rightarrow X_{t} \geq \beta(t).
\end{align*}
Thus, on the interval $[0,\tau_{\beta}]$, $X_{t}\in[\inf_{t\geq 0}\beta(t),\infty)$. Let $b=\inf_{t\geq 0}\beta(t)$. This implies
\begin{align*}
    \int_{0}^{\tau_{\beta}}\left(\frac{\sigma}{2}-\frac{I}{\sigma}-\frac{V_{r}}{\tau \sigma}+\frac{1}{\tau \sigma}e^{-\sigma X_{s}}\right)^{2}ds \leq \left(\frac{\sigma}{2}-\frac{I}{\sigma}-\frac{V_{r}}{\tau \sigma}\right)^{2}\tau_{\beta} +\frac{1}{\tau^{2}\sigma^{2}}e^{-2\sigma b}\tau_{\beta}+\frac{2}{\tau \sigma}e^{-\sigma b}\tau_{\beta}<\infty\,.
\end{align*}
Since the threshold $\beta(t)$ lies below the process $X$, the Exact simulation method requires an update. We simulate the first-passage time of the Brownian motion using the algorithm described in Appendix B and apply the Exact Algorithm with the updated Bessel process:
\[
R_s = X_{\tau_{\beta}-s}^{\mathbb{Q}}-\beta(\tau_{\beta}-s).
\]
This updated approach is used to generate spike trains for varying input currents, as shown in Figure \ref{fig:enter-label}.
We observe that with a higher input current, there is a higher number of spikes. Also, when $I=0$, the spike time position is very random, but as $I$ increases, the spike time position becomes more similar. For instance, when $I=20$, the spike time position is very similar between the times 0 s and 0.5 s.\HLZ\\
Notably, in the case of a quadratic SDE, we can use our Exact method. This is a substantial advantage of the Exact method, as other time-discretization methods like the Euler-Maruyama method cause rapid changes in the value of process $V$, which leads to instabilities, cf. \cite{szpruch2010numerical}. For example, when $V$ becomes large, the term $-\frac{1}{\tau}V_t(V_t-V_r)$ grows quickly, and the Euler-Maruyama method will struggle to handle this large change accurately. Even small errors in the drift could accumulate, potentially causing blow-ups or oscillations in the solution.
\section{Conclusion}
In the present work, we extended an Exact simulation method for the first-passage time of stochastic differential equations to time-dependent thresholds. The Exact simulation method generates the FPT of the process by accepting or rejecting samples of the FPT of a Brownian motion. We derived the rejection probability for sampling  FPT to a time-dependent threshold, using Girsanov's transformation. When the explicit density of the FPT of the Brownian motion is available, we sample directly from it and use the rejection probability to accept or reject the sample for the FPT of our SDE. If the explicit density is unavailable, we adopt the approximation technique described in \cite{herrmann2016first}. Furthermore, we studied the time complexity of the Exact Algorithm and observed an exponential bound in the number of iterations. To address this, we proposed techniques that help to reduce the computational time of the algorithm and also proved that with this technique, the bound in the number of iterations becomes linear.\HLZ\\
One limitation of our method is the complexity of directly evaluating the rejection probability, as it involves an integral of a stochastic process with respect to time. To overcome this, we applied a technique that bypasses the need to compute the probability explicitly by determining a rejection region and using the thinning method (cf. \cite{lewis1979simulation}) to perform acceptance-rejection. If all points lie outside the rejection region, the sample is accepted. 
A further limitation is that extending the method to multidimensional settings is possible, but it would require assumptions that are rarely practical, cf. \cite{10.1007/978-3-642-27440-4_7}.\HLZ\\
Our approach is particularly well-suited for applications, where accurate first-passage time predictions are crucial.\HLZ\\
Future work in this direction could be extending the method to handle jump processes and scenarios with two-sided, time-dependent thresholds. 
\section*{Acknowledgements}
Devika Khurana thanks the Austrian Science Fund (FWF) for financial support through the project "Stochastic Models and Methods for the Study of Olfaction" - Grant-DOI 10.55776/I4620.
\section*{Appendix A. Bessel Bridge}
In this appendix, we explain the Bessel Bridge and describe the Proposition that proves that the process $R$ that we used in Section 2.2 is a Bessel Bridge.\HLZ\\ 
\textbf{Definition A.1} \textit{Bessel Process $R=\{R_t,t\geq0\}$}:
\begin{itemize}
    \item[1.] Let $W = (W_1,W_2,...,W_n)$ be a n-dimensional Brownian motion in $R^n$, $n\geq 2$. Then the distance $R_t$ of $W$ from the origin i.e.
\[
R_t=\sqrt{W_{1}(t)^{2} + W_{2}(t)^{2} + ... + W_{n}(t)^{2}}
\]
is known as the n-dimensional Bessel process.
\item[2.] $R$ is the solution to the following SDE:
\[
dR_t = \frac{n-1}{2R_t}dt + dB_t,
\]
where $B_t$ is a n-dimensional Brownian motion.
\end{itemize}
\textbf{Definition A.2} \textit{Bessel Bridge $BR=\{BR_t,t\geq0\}$}: A Bessel Bridge $BR$ is a Bessel process that is conditioned to start at a specific value at $t=0$ and to reach a predetermined value at a later time $T>0$.\HLZ\\
For a more detailed overview of these processes, please refer to \cite{oksendal2013stochastic,revuz2013continuous}.\HLZ\\
\textbf{Proposition A.1} The process $R_t = \beta(\tau_{\beta}-t)-W_{\tau_{\beta}-t}$; $t\in[0,\tau_{\beta}]$ defined in Section 2.2, where $W$ is a 1-D Brownian motion starting at $x_0$ and ends at $\beta(\tau_{\beta})$ is a Bessel Bridge.\\

\noindent For a detailed overview and proof, please refer to \cite{hernandez2013hitting}.

\section*{Appendix B. Algorithm to generate $\tau_{\beta}^{W}$, when $\beta(t)$ is curvy.}
In what follows, we give the description of the algorithmic approach to approximate the first-passage time of Brownian motion to a curvy threshold $\beta(t)$ as introduced in \cite{herrmann2016first}.

\begin{table}[ht]
    \begin{tabular}{l}
    \hline
        Generating $\tau_{\beta}^{W}$ for a curvy threshold $\beta(t)$, $\beta(0)>0$\\
    \hline
         Let, $\epsilon>0$ be a small parameter and $r$ be any real number, such that $r>\rho$,\\ where $\rho \leq \inf_{t\geq0} \beta'\left(t\right)$.\\
    Initialization: $T=0$, $H=\beta\left(0\right)$\\
    \\
    While $H>\epsilon$ and $T<K$\\
    \textbf{Step 1} Generate an inverse Gaussian random variable, $\hat{G}$ with parameters\\ \hspace{1.3cm}$H/r$ and $H^2$.\\
    \textbf{Step 2} Update the values\\
       \hspace{1.3cm}$H\leftarrow \beta(T+\hat{G}) - \beta\left(T\right) + r\hat{G},$\\
       \hspace{1.3cm}$T\leftarrow\hat{G}+T$\\
    \textbf{Outcome}: A sample of  $\tau_{\beta}^{W}\leftarrow T \wedge K$\\
    \hline
  \end{tabular}    
\end{table}

If $\beta(0)<0$, meaning the threshold starts below the Brownian path and $\beta(t)$ is non-decreasing, we employ tilted horizontal lines directed towards the path. The algorithm in this case is similar to the previous approach, with new $r$ and $H$:
\[
r>\rho=\sup_{t\geq 0 }\beta'(t).
\]
We initialize $H$ as $\beta(0)$ and update it as follows in Step 2 of the algorithm:
\[
H\leftarrow r\hat{G} + \beta\left(T\right) -\beta(T+\hat{G}) 
\]

\bibliographystyle{plain}
\bibliography{Ref}
\end{document}